\documentclass[12pt]{article}

\usepackage{latexsym}
\usepackage{setspace}
\usepackage{color}
\usepackage{float}

\usepackage{amsmath}
\usepackage{amssymb}
\usepackage{latexsym}
\usepackage{amsmath}
\usepackage{mathrsfs}
\usepackage{amssymb}
\usepackage{natbib}
\usepackage{multirow}
\usepackage{subfigure}
\usepackage{subcaption}
\usepackage{epstopdf}
\usepackage{pgfplots} 
\usepackage{subfigure}
\setcitestyle{authoryear,open={(},close={)}}
\usepackage{hyperref}
\hypersetup{colorlinks,citecolor=blue}
\usepackage{graphics}
\usepackage{graphicx}
\makeatletter
\def\ps@pprintTitle{%
  \let\@oddhead\@empty
  \let\@evenhead\@empty
  \let\@oddfoot\@empty
  \let\@evenfoot\@oddfoot
}
\makeatother
\makeatletter \oddsidemargin  -.5cm \evensidemargin -.5cm \textwidth
17.5cm \topmargin -.65in \textheight 24cm

\newcommand{\bd}{\begin{definition}}
\newcommand{\ed}{\end{definition}}
\newcommand{\br}{\begin{remark}}
\newcommand{\er}{\end{remark}}
\newcommand{\bea}{\begin{eqnarray}}
\newcommand{\eea}{\end{eqnarray}}
\newcommand{\beann}{\begin{eqnarray*}}
\newcommand{\eeann}{\end{eqnarray*}}
\newtheorem{theorem}{Theorem}[section]

\newtheorem{lemma}[theorem]{Lemma}

\newtheorem{corollary}[theorem]{Corollary}
\newtheorem{remark}{Remark}[section]

\numberwithin{equation}{section}
\numberwithin{equation}{section}
\title{Estimating ordered variance of two scale mixture of normal distributions}
\author{Shrajal Bajpai and Lakshmi Kanta Patra\footnote
	{\baselineskip=10pt
		~lkpatra@iitbhilai.ac.in; patralakshmi@gmail.com}  \\
	\it  Department of Mathematics\\
	\it Indian Institute of Technology Bhilai, Durg, India-491002}
\begin{document}
	\date{}
	\maketitle
	\begin{abstract}
		 This study investigates component wise estimation of ordered variances of scale mixture of two normal distributions. For this study two special loss functions are considered namely squared error loss function and entropy loss function. We have derived the general improvement results  and based on these results the estimators that outperform BAEE are obtained. Moreover under certain sufficient conditions a class of improved estimators is proposed for both loss functions. As a special case of scale mixture of normal distribution the results are applied to the multivariate t-distribution and obtained the improvement results. For this case a detailed numerical comparison is carried out which validates our theoretical findings.

	\end{abstract}

\section{Introduction}
 In various real-life scenarios, the experimenter often confronts prior information about the restrictions on parameters. For example, humidity levels in cities located near coastal regions are generally higher than those in inland towns. It is natural to assume the average pollution level of an industrial city is higher than that of a non-industrial city. Thus, it is natural to take some order restrictions on the parameter in statistical models. Numerous researchers studied the problem of estimating order-restricted parameters of various probability distributions. Some earlier works in this direction are on developing a maximum likelihood estimator when the parameters satisfy an order restriction. For a relevant discussion in this direction, we refer to  \cite{barlow1972statistical}, \cite{wright1988order}, \cite{shi1998maximum}, \cite{marchand2004estimation} and \cite{van2006minimaxity}. 

In the past, several authors discussed the problem of finding improved estimators of the ordered parameters of various probability distributions. For some initial work on finding improved estimators of the ordered parameters of normal and exponential distributions, we refer to \cite{kumar1992inadmissibility}, \cite{misra1994estimation},        \cite{pal1988decision}, \cite{vijayasree1995componentwise}, \cite{misra1997estimation}, \cite{misra2002smooth}, \cite{oono2005estimation}, \cite{kumar2005james}.

\cite{kumar2008estimating} investigated the problem of estimating the restricted normal means. The authors proposed several estimators and established the inadmissibility of the usual estimator. \cite{tsukuma2008stein}  discussed the estimation of restricted normal mean vector with known variance under the quadratic loss function when means are restricted to a polyhedral convex cone. They obtained a generalized Bayes estimator and studied its properties. The authors proved that the generalized Bayes estimator is minimax and admissible in the one-or two-dimensional case, but can be improved by a shrinkage estimator when the dimension is greater than two. \cite{tripathy2013estimating} discussed the estimation of the common variance of two normal populations with ordered means under the scale invariant quadratic loss function. They proved the inadmissibility of usual estimators by proposing some dominating estimators. 

Estimation of ordered normal means under Linex loss has been studied by \cite{ma2013estimation}. They proved that under certain sufficient conditions, the plug-in estimator dominates the unrestricted MLE. \cite{jena2021alternative}  considered the problem of estimating the common mean of two normal populations having ordered variance. They proposed estimators that dominate the usual estimator. The estimation of ordered means of a bivariate normal distribution has been studied by \cite{patra2017estimating} under a sum-squared error loss function. The authors considered the estimation problem for both cases, that is, equal and unequal variances. For each case, a class of minimax estimators is proposed that improves upon the usual estimators. \cite{patra2021minimax} considered the minimax estimation of common variance and common precision of two normal populations with ordered means under an invariant loss function. In this work, the authors showed the inadmissibility of the best affine equivariant estimator and MLE by proposing various dominating estimators. \cite{jana2022estimation} investigated the estimation of the ordered variances of two normal populations with a common mean under the quadratic loss function. BAEE and its improvements are derived for ordered variances. Moreover, the authors estimated the smallest standard deviation of $k(\ge 2)$ normal populations. They demonstrated the inadmissibility of   BAEE for the smallest standard deviation by proposing improved equivariant estimators. \cite{jena2023point} explored the point and interval estimation of powers of scale parameters of two normal populations having a common mean. Maximum likelihood and plug-in estimators have been proposed. The authors have proved inadmissibility of the best affine equivariant estimator. Moreover, in this work, the authors have proposed several interval estimators.

In recent years, the improved estimation of the scale parameter in mixture models has been studied by \cite{petropoulos2005estimation}, \cite{petropoulos2010class}, \cite{petropoulos2017estimation}, and \cite{patra2021improved}. 
The ordered parameter estimation for the mixture of normal distribution models has not been previously studied.  In this article, we consider estimating ordered variances of scale mixture of multivariate normal distribution.
%
Let $Y_1,\dots, Y_n$ are $p$ dimensional random vector having  joint density 
\begin{eqnarray}\label{model1}
	f(y_1,\dots,y_n,\theta,\sigma^2)=\int \frac{\tau^{np/2}}{(2\pi)^{np/2}\sigma^{np}}\exp \left(-\frac{\tau}{2\sigma^2}\sum_{i=1}^{n}(y_i-\theta)^t(y_i-\theta)\right)d\Pi(\tau)
\end{eqnarray}
Interpretation of (\ref{model1}), is for given $\tau>0$ we have $Y_1,\dots, Y_n$ are i.i.d with normal $N_{p}(\theta,\frac{\sigma^2}{\tau}I_p)$. This
model is named as scale mixture of normal distribution. Here $\tau$ is the mixing parameter having distribution $\Pi(\tau)$. So for a given $\tau$, sufficiency reduces to $\bar{Y} \sim N_p(\theta,\frac{\sigma^2}{n\tau}I_p)$,  $S \sim \frac{\sigma^2}{\tau} \chi^2_{p(n-1)}$  and both are independent.

In this article we consider two scale mixture of normal models as, given $\tau>0$ 
\begin{equation}\label{model2}
	\boldsymbol{X_1}\sim N_p(\boldsymbol{\theta}_1,\frac{\sigma_1^2}{\tau} I) , ~~S_1 \sim\frac{ \sigma_1^2}{\tau}\chi ^2_{m_1}~~\mbox{ and }~~\boldsymbol{X_2}\sim N_q(\boldsymbol{\theta}_2,\frac{\sigma_2^2}{\tau} I ),~~S_2 \sim \frac{\sigma_2^2}{\tau}\chi ^2_{m_2}
\end{equation}
with the order restriction $\sigma_1 \le \sigma_2$ and $\boldsymbol{\theta}_1 \in \mathbb{R}^p,\boldsymbol{\theta}_2 \in \mathbb{R}^q$ and $\boldsymbol{X_1}$, $S_1$, $\boldsymbol{X_2}$ and $S_2$ are independently distributed.  We take $T_1=\|\boldsymbol{X_1}\|^2$ and $T_2=\|\boldsymbol{X_2}\|^2$. Then 
$$T_1 \sim \frac{\sigma_1^2}{\tau}\chi^2_p(\nu_1), ~\nu_1=\frac{\tau\|\boldsymbol{\theta_1}\|^2}{\sigma^2_1}~~\mbox{ and }~~T_2\sim  \frac{\sigma^2_2}{\tau}\chi^2_q(\nu_2), ~\nu_2=\frac{\tau\|\boldsymbol{\theta_2}\|^2}{\sigma^2_2}.$$
So $T_1|L=l \sim \frac{\sigma_1^2}{\tau}\chi^2_{p+2l}$ and $T_2|K=k \sim \frac{\sigma_2^2}{\tau}\chi^2_{q+2k}$, where $L$ and $K$ are Poisson random variable with parameter $\nu_1$ and $\nu_2$ respectively. We will consider component wise estimation of $\sigma_i^2$, $i=1,2$ with respect to the following loss functions. 
\begin{align*}
&\mbox{Squared error loss~}: ~L_1(\delta,\sigma^2)=\left(\delta/\sigma^2-1\right)^2 \\
&\mbox{Entropy loss~}: ~L_2(\delta,\sigma^2)= (\delta/\sigma^2) - \ln(\delta/\sigma^2)-1
\end{align*}
At first we use the principal of invariance. So we consider the affine group of transformations $$\mathcal{G}_{a_i,b_i}=\left\{g_{a_i,b_i}(x)=a_ix+b_i,i=1,2, a_i>0, b_i \in \mathbb{R}\right\}$$ 
Under this group, the form of the affine equivariant estimators of $\sigma_i^2$ can be obtained as 
\begin{equation}
	\delta_{\alpha_{i}}(\boldsymbol{X_i},\underline{S_i})=\alpha_{i} S_i
\end{equation}
where $\alpha_i>0$ is a constant. The following lemma provides the BAEE of $\sigma^2_i$.
\begin{lemma}
	\begin{enumerate}
		\item[(i)] For the squared error loss function $L_1(t)=(t-1)^2$, the BAEE of $\sigma_i^2$ is $\delta^Q_{0i}=\alpha_{iQ}S_i$, where $\alpha_{iQ} =\frac{E(\tau^{-1})}{(2+m_i)E(\tau^{-2})}$.
		\item [(ii)] For the entropy loss function $L_2(t)=t-\ln t-1$, the BAEE of $\sigma_i^2$ is $\delta^E_{0i}= \alpha_{iE}S_i$ where $\alpha_{iE}= \frac{m_i}{E(\tau^{-1})} .$
	\end{enumerate}
\end{lemma}

The rest of the paper is organized as follows. In Section \ref{sec2}, we have discussed the improved estimation of $\sigma_1^2$. We have proposed estimators that dominate the BAEE under the squared error loss function and the entropy loss function. Furthermore, a class of estimators that improves upon the BAEE is given. In Section \ref{sec3}, estimation of $\sigma_2^2$ is considered. Improvement results similar to those in Section \ref{sec2} are also obtained for this case. As an application we have obtained improved estimators of the ordered scale parameters of multivariate t- distribution. Finally, a simulation study has been carried out to asses the risk performance of the proposed estimators.

\section{Estimation of $\sigma_1^2$} \label{sec2}
In this section, we will consider the improved estimation of $\sigma^2_1$. Here, we are going to prove that the BAEE of $\sigma^2_1$ is inadmissible by deriving an improved estimator. Consider a broader class of estimators based on $S_1$ and $S_2$. 
$$\mathcal{C}_1 = \left\lbrace {\delta_{\phi_{1}}= \phi_{1}(Z_1)S_1 \big|  Z_1= \frac{S_2}{S_1}},\   \phi_1\  \text {is positive measurable function}\right\rbrace. $$
It is observed that for BAEE $\delta_{01}^Q$ and $\delta_{01}^E$, $\phi_{1}(Z_1)= \alpha_{1Q},$ and $\alpha_{1E}$ respectively.
Risk function of the estimator $\delta_{\phi_{1}}$  under squared error loss function is 
\[
R(\delta_{\phi_{1}},\sigma^2_1,\sigma^2_2) = E\Big[\Big(\frac{\phi_1(Z_1)S_1}{\sigma_1^{2}}-1\Big)^{2}\Big]
=\; E\left[ E\left\lbrace \Big(\frac{\phi_1(Z_1)S_1}{\sigma_1^{2}}-1\Big)^{2}\Big\vert Z_1\right\rbrace \right] .
\]
To find a minimizer of the risk of $\delta_{\phi_1}$, we will find the minimzer of the inner conditional risk. For a given $\tau>0$, the joint density of $V_1$ and $Z_1$ is proportional to 
$$\eta^{\frac{m_2}{2}+1} {v_1^{\frac{m_1+m_2}{2}-1}}{z_1^{-(2+\frac{m_1}{2}-1)}} e^{-v_1(1+\eta^{-2}z_1)},$$
where $\eta=\frac{\sigma_1}{\sigma_2} \le 1$. 
After some simplification we get the minimizer as 
\begin{align}\label{e1}
	\phi_1^Q(z_1,\eta)& =\frac{\dfrac{1}{\sigma_1^{2}}\; E\big(S_1 \mid Z_1=z_1\big)}
	{\dfrac{1}{\sigma_1^{4}} E\big(S_1^{2} \mid Z_1=z_1\big)}
	= \frac{%
		E\!\left[\tau^{-1}\,E\!\left(\left.\frac{\tau}{\sigma_1^2}\,S_1\ \right|\ Z_1=z_1\right)\right]%
	}{%
		E\!\left[\tau^{-2}\,E\!\left(\left.\frac{\tau^2}{\sigma_1^4}\,S_1^2\ \right|\ Z_1=z_1\right)\right]%
	}.
\end{align}
Denote  $V_1=\frac{\tau}{\sigma_1^2}\,S_1$, then the equation (\ref{e1}) can be written as
$$	\phi_1^Q(z_1,\eta) =\frac{
	E\!\left[\tau^{-1}\,E\!\left(V_1 | Z_1=z_1\right)\right]
}{E\!\left[\tau^{-2}\,E\!\left(V_1^2 | Z_1=z_1\right)\right]
}.$$
After simplifying we get 
$$	\phi_1^Q(z_1,\eta)=\frac{(1+\eta^{-2} z_1)E(\tau^{-1})}{({m_1+m_2}+2)E(\tau^{-2})}\le \frac{(1+ z_1)E(\tau^{-1})}{({m_1+m_2}+2)E(\tau^{-2})}=\phi_{11}^{Q}(z_1)~\mbox{ (say)}.$$
The following theorem gives an improvement result over the estimator $\delta_{\phi_{1}}$.
\begin{theorem}\label{th1}
	Under the squared error loss function $L_1(\cdot)$, the risk of the estimator
	\begin{equation*}
		\delta^Q_{S1}(X,S)=\min \left\{\phi_1(Z_1),\phi_{11}^{Q}(Z_1)\right\}S_1,
	\end{equation*}
	is nowhere larger than that of the estimator $\delta_{\phi_1}$ provided $P\left(\phi^Q_{11}(Z_1)<\phi_1(Z_1)\right)>0$.
\end{theorem}

\noindent Similar to the squared error loss function $L_1(.)$, the minimzer of the risk function of $\delta_{\phi_1}$ under entropy loss function $L_2(.)$  is obtained as 
\[\phi_1^{E}(z_1) =\frac{1+\eta^{-2}z_1}{E\left[ {\tau^{-1}}(m_1+m_2)\right] } \le \frac{1+z_1}{E\left[ {\tau^{-1}}(m_1+m_2)\right]}  = \phi_{11}^{E}(z_1)~~\mbox{ (say)}.\]
\begin{theorem}\label{th2}
	The estimator
	\begin{equation*}
		\delta^E_{S1}(X,S)=\min \left\{\phi_1(Z_1),\phi_{11}^{E}(Z_1)\right\}S_1,
	\end{equation*}
	dominates $\delta_{\phi_{1}}$ with respect to entropy  loss function $L_2(\cdot)$ provided $P\left(\phi^E_{11}(Z_1)<\phi_1(Z_1)\right)>0$.
\end{theorem}
As an application of Theorems \ref{th1} and \ref{th2}, in the following corollary we will give estimator that dominates BAEE. 
\begin{corollary}
	\begin{itemize}
		\item [(i)] For the squared error loss function $L_1(.)$ the estimator 
		\begin{equation}
			\delta_{11}^Q = 
			\begin{cases}
				\frac{(1+ Z_1)E(\tau^{-1})}{({m_1+m_2}+2)E(\tau^{-2})}S_1, 
				& \text{if } 0 \le Z_1 \le \dfrac{m_2}{m_1 + 2}, \\[1.2em]
				\frac{E(\tau^{-1})S_{1}}{(2+m_1)E(\tau^{-2})}, 
				& \text{otherwise,}
			\end{cases}
		\end{equation}
		dominates $\delta_{01}^Q$.
		\item [(ii)] 	Under the entropy loss function $L_2(.)$ the estimator 
		\begin{equation}
			\delta_{11}^E= 
			\begin{cases}
				\dfrac{(1+Z_1)\,S_1}{(m_1 + m_2)E(\tau^{-1})}, 
				& \text{if } 0 \le Z_1 \le {\frac{m_2}{m_1}}, \\[1.2em]
				\frac{S_1}{m_1E(\tau^{-1})}, 
				& \text{otherwise,}
			\end{cases}
		\end{equation}
		has uniformly smaller risk than $\delta_{01}^E$.
	\end{itemize}
\end{corollary}
\noindent Next we explore a large class of estimators utilizing the information contained in statistics $S_1, S_2$ and $\boldsymbol{X}_1$. Consider a class of estimators as
$$\mathcal{C}_2= \left\lbrace \delta_{\phi_{2}}= \phi_{2}(Z_1,Z_2)S_1 \  \big|  \ Z_1=\frac{S_2}{S_1},Z_2 =\frac{||\boldsymbol{X}_1||^2}{S_1}, \phi_{2}  \ \text{is positive measurable function}\right\rbrace.$$
Under the squared error loss function $L_1(.)$,  risk of the $\delta_{\phi_{2}}$ is minimized at
\begin{equation}\label{re1}
	\phi_{2}^Q(\eta,z_1,z_2) 	= \frac{%
		E \left[\tau^{-1}\,E\!\left(\left.\frac{\tau}{\sigma_1^2}\,S_1 \right|\ Z_1=z_1,  Z_2=z_2, L=l \right)\right]%
	}{%
		E \left[\tau^{-2}\,E\!\left(\left.\frac{\tau^2}{\sigma_1^4}\,S_1^2 \right|\ Z_1=z_1, Z_2=z_2, L=l\right)\right]%
	}.
\end{equation}
Where $L$ is the Poisson distribution with parameter $\nu_1$, for given $\tau>0$, and $L=l$  the joint probability density function of $V_1$ given  $Z_1$ and $Z_2$  is $\propto$
$v_1^{\frac{m_1+m_2+p+l}{2}-3}e^{-\frac{v_1}{2}(1+\eta^{-2}z_1+z_2)}$, $\eta \le 1$ and $v_1>0$. Consequently we obtain 
\begin{equation}\label{re2}
	\phi_{2}^Q(\eta,z_1,z_2) 	= \frac{%
		E \left[\tau^{-1}\int_{0}^{\infty}v_1^{\frac{m_1+m_2+p+2l}{2}}e^{-\frac{v_1}{2}(1+\eta^{-2}z_1+z_2)}dv_1\right]%
	}{%
		E \left[\tau^{-2}\int_{0}^{\infty}v_1^{\frac{m_1+m_2+p+2l}{2}+1}e^{-\frac{v_1}{2}(1+\eta^{-2}z_1+z_2)}dv_1\right]%
	}.
\end{equation}
After simplifying equation (\ref{re2}), we get
	
	\begin{align*}
		\phi_{2}^Q(\eta,z_1,z_2) &=\frac{( 1+\eta^{-2}z_1+z_2)}{m_1+m_2+p+2l+2} \frac{E(\tau^{-1})}{E(\tau^{-2})}\\&\le \frac{( 1+z_1+z_2)}{(m_1+m_2+p+2l+2)}  \frac{E(\tau^{-1})}{E(\tau^{-2})}	
		\le \frac{( 1+z_1+z_2)}{(m_1+m_2+p+2)} \frac{E(\tau^{-1})}{E(\tau^{-2})}\\
		&=\phi_{21}^Q(z_1,z_2) ~~~\text{(say)}.
	\end{align*}

	
	\begin{theorem}\label{Th3}
		Under the squared error loss function $L_1(\cdot)$, the risk of the estimator
		\begin{equation*}
			\delta^Q_{S2}=\min \left\{\phi_2(Z_1,Z_2),\phi_{21}^{Q}(Z_1,Z_2)\right\}S_1,
		\end{equation*}
		is nowhere larger than that of the estimator $\delta_{\phi_2}$ provided $P\left(\phi^Q_{21}(Z_1,Z_2)<\phi_2(Z_1,Z_2)\right)>0$.
	\end{theorem}
	
	\noindent Again under the entropy loss function minimizer of the risk function  is   obtained as 
	\begin{equation*}
		\phi_{2}^{E}=\left[ E\left({\tau^{-1}}\;E\big(V_1 \mid Z_1=z_1,\; Z_2=z_2,\;  L=l\big)\right)\right]^{-1}.
	\end{equation*}
	Consequently, we have
	\begin{align*}
		\phi_{2}^{E}& = \frac{(1+\eta^{-2} z_1+z_2)}{m_1+m_2+p+2l-4}\left[  E\left( \tau^{-1}\right) \right] ^{-1}
		\le \frac{(1+ z_1+z_2)}{(m_1+m_2+p+2l-4)}\left[  E\left( \tau^{-1}\right) \right] ^{-1}\\
		&\le \frac{(1+z_1+z_2)}{(m_1+m_2+p-4)}\left[  E\left( \tau^{-1}\right) \right] ^{-1} = \phi_{21}^{E}(Z_1,Z_2) \ \ \text {(say)}.
	\end{align*}
	
	\begin{theorem}\label{Th4}
		The risk of the estimator
		\begin{equation*}
			\delta^E_{S2}=\min \left\{\phi_2(Z_1,Z_2),\phi_{21}^{E}(Z_1,Z_2)\right\}S_1,
		\end{equation*}
		is uniformly smaller than that of the estimator $\delta_{\phi_2}$ with respect to  entropy loss function $L_2(\cdot)$, provided $P\left(\phi^E_{21}(Z_1,Z_2)<\phi_2(Z_1,Z_2)\right)>0$.
	\end{theorem}

	
	\begin{corollary}
		\begin{itemize}
			\item [(i)] For the squared error loss function $L_1(.)$ the estimator 
			\[\delta_{12}^Q =
			\begin{cases}
				\frac{1 + Z_1 + Z_2}{\,m_1 + m_2 + p + 2\,}
				\;\frac{E(\tau^{-1})}{E(\tau^{-2})},
				& \text{if } 
				0 \le Z_2 \le \dfrac{m_2 + p}{m_1 + 2} - Z_1
				\text{ and }
				0 \le Z_1 \le \dfrac{m_2}{m_1 + 2},
				\\[1.2em]
				
				\frac{S_1}{m_1 + 2}\;
				\frac{E(\tau^{-1})}{E(\tau^{-2})},
				& \text{otherwise},
			\end{cases}
			\]
			dominates $\delta_{01}^Q$.
			\item [(ii)] The estimator 
			\[
			\delta_{12}^E =
			\begin{cases}
				
				\frac{1 + Z_1 + Z_2}{\,m_1 + m_2 + p - 4\,}
				\;\frac{1}{E(\tau^{-1})},
				& \text{if } 
				0 \le Z_2 \le \dfrac{m_2 + p-4}{m_1} - Z_1
				\text{ and }
				0 \le Z_1 \le \! \frac{m_2}{m_1},
				\\[1.2em]
				
				\frac{S_1}{m_1}\;
				\frac{1}{E(\tau^{-1})},
				& \text{otherwise},
			\end{cases}
			\]
			dominates  $\delta_{01}^E$ with respect to entropy loss function $L_2(.)$.
		\end{itemize}
	\end{corollary}
	\subsection{A class of improved estimators}\label{sec2.1}
	In this section, we will obtain a class of improved estimators using the approach proposed by \cite{kubokawa1994unified}.   
	For a given $\tau>0$ we denote $g_i$ and $G_i$ be the density function and distribution function of $V_i$, $i=1,2$.
	The following theorem provides certain sufficient conditions 
	under which the class of estimators $\delta_{\phi_{1}}= {\phi_1(Z_1)}S_1$ improves upon the BAEE under squared error loss function.
	\begin{theorem}\label{kth1}
		The risk of estimator $\delta_{\phi_{1}}= {\phi_1(Z_1)}S_1$ is nowhere larger than $\delta_{01}^Q$ under squared error loss function $L_1(.)$ provided $\phi_{1}(u_1)$ satisfies the following conditions:
		\begin{enumerate}
			\item[(i)]$\phi_1(u_1)$ is non-decreasing and $\lim_{u_1\rightarrow\infty}\phi_1(u_1) = {\alpha_{1Q}},$
			\item [(ii)]$\phi(u_1) \ge \phi^*_1(u_1) = \frac{1}{m_1+m_2+2}\frac{E[\tau^{-1}]\displaystyle \int_{0}^{u_1}
				t^{\frac{m_2}{2}-1}(1+t)^{-(\frac{m_1+m_2}{2}+1)}dt  }{E[\tau^{-2}]\displaystyle \int_{0}^{u_1} t^{\frac{m_2}{2}-1}(1+t)^{-(\frac{m_1+m_2}{2}+2)}dt \ }.$
			
		\end{enumerate}
	\end{theorem}
	\noindent 	\textbf{Proof:} The idea of proof of this theorem is similar to \cite{petropoulos2017estimation}. So we omit it.
		\begin{theorem}\label{kth2}
			The estimator $\delta_{\phi_{1}}= {\phi_1(Z_1)}S_1$ dominates the estimator $\delta_{01}^{E}$ with respect to entropy loss function $L_2(.)$ provided $\phi_{1}(u_1)$ satisfies the following conditions:
			\begin{enumerate}
				\item[(i)]$\phi_1(u_1)$ is non-decreasing and $\lim_{u_1\rightarrow\infty}\phi_1(u_1) = {\alpha_{1E}},$
				\item [(ii)]$\phi_1(u_1) \ge \phi^*_1(u_1) = \frac{1}{m_1+m_2}\frac{\displaystyle \int_{0}^{u_1}
					t^{\frac{m_2}{2}-1}(1+t)^{-(\frac{m_1+m_2}{2})}dt  }{E[\tau^{-1}]\displaystyle \int_{0}^{u_1} t^{\frac{m_2}{2}-1}(1+t)^{-(\frac{m_1+m_2}{2}+1)}dt \ }.$
				
			\end{enumerate}
		\end{theorem}
		\noindent 	\textbf{Proof:}
		The risk difference of the estimators  $\delta_{01}^E$ and $\delta_{\phi_{1}}= {\phi_1(Z_1)}S_1$  is 
		\begin{align*}
			RD
			& =\; E \left[\,\frac{\alpha_{1E} S_1}{\sigma_1^2} - \ln\!\left(\frac{\alpha_{1E} S_1}{\sigma_1^2}\right) - 1 \,\right]-E  \left[\,\frac{\phi_1(Z_1) S_1}{\sigma_1^2} - \ln\!\left(\frac{\phi_1(Z_1) S_1}{\sigma_1^2}\right) - 1 \,\right]\\
			&= E  \left[\,\alpha_{1E} Y_1 - \ln(\alpha_{1E} Y_1) - 1 \,\right]-E\left[\,\phi_1(Z_1) Y_1 - \ln(\phi_1(Z_1) Y_1) - 1 \,\right]
		\end{align*}
		which  can be written as 
		\begin{align*}
			%
			RD	&= E\Bigg[\int_{1}^{\infty}
			\phi_1'\Big(\eta^{-2}t\tfrac{Y_{2}}{Y_{1}}\Big)\Big(\eta^{-2}Y_{2}\Big)
			-\frac{\phi_1'\Big(\eta^{-2}t\tfrac{Y_{2}}{Y_{1}}\Big)\big(\eta^{-2}Y_{2}\big)}
			{\phi_1\Big(\eta^{-2}t\tfrac{Y_{2}}{Y_{1}}\Big)Y_{1}} \,dt\Bigg]\\
			&=\eta^{-2}\int_{0}^{\infty}\!\int_{0}^{\infty}\!\int_{0}^{\infty}\!\int_{1}^{\infty}\!
			\phi_1'\!\Big(\frac{\eta^{-2} t y_{2}}{y_{1}}\Big)
			\Big(1-\frac{1}{\phi_1\!\big(\frac{\eta^{-2}ty_{2}}{y_1}\big)\,y_{1}}\Big)y_2
			\tau g_{1}(\tau y_{1})\,\tau g_{2}(\tau y_{2})\\
			&~~~~~~~~~~~~~~~~~~~~~~~~~~~~~~~~~~~~\times\,dt\,dy_{2}\,dy_{1}\,d\Pi(\tau).
		\end{align*}
		Making the transformations \(u_{1}=\dfrac{t y_{2}}{y_{1}}\) and $\frac{u_1y_1}{t}=x$ then we obtain 
		\begin{align*}
			RD&\ge \int_{0}^{\infty}\!\int_{0}^{\infty}\!\int_{0}^{\infty}\int_{1}^{\infty}
			\frac{u_1y_1^2}{t^2}\,\phi_1'\big(\eta^{-2}u_{1}\big)
			\Big(1-\frac{1}{\phi_1\big(\eta^{-2}u_{1}\big)\,y_{1}}\Big)
			\tau g_{1}(\tau y_{1})\,\tau g_{2}(\tau \frac{u_1y_1}{t})\\
			&~~~~~~~~~~~~~~~~~~~~~~~~~~~~~~~~~~~\times dt\,du_{1}\,dy_{1}\,d\Pi(\tau).\\
		%
		%
			&\ge\int_{0}^{\infty}\int_{0}^{\infty}\int_{0}^{\infty}\int_{0}^{u_1y_1}
			x\,\phi_1'(\eta^{-2}u_1)\Bigg(1-\frac{1}{\phi_1(\eta^{-2}u_1)\,y_{1}}\Bigg)
			\tau\,g_{1}(\tau y_{1})\tau g_{2}(\tau x)\,
			\frac{t}{u_1}\;\\
			&~~~~~~~~~~~~~~~~~~~~~~~~~~~~~~~~~~~ \times dx\,du_{1}\,dy_{1}\,d\Pi(\tau)\\
			&\ge\int_{0}^{\infty}\phi_1'(\eta^{-2}u_1)\int_{0}^{\infty}\int_{0}^{\infty}
			\left(1-\frac{1}{\phi_1(\eta^{-2}u_1)\,y_{1}}\right) y_{1}\,\tau\,g_{1}(\tau y_{1}) \int_{0}^{\tau u_1y_1}g_{2}(x)\\
			&~~~~~~~~~~~~~~~~~~~~~~~~~~~~~~~~~~~\times \;dx  \ du_1\,dy_{1} d\Pi(\tau).
		\end{align*}
		So to ensure $RD\ge 0$ we must have
		\[
		\int_0^\infty\int_0^\infty \Big( y_1 - \frac{1}{\phi_1(\eta^{-2} u_1)}\Big)\,\tau g_1(\tau y_1)\,G_2(\tau u_1y_1)\,dy_1\,d\Pi(\tau)\ge 0,
		\]
		%
		%

		that is 
		
		\begin{equation*}
			\phi(u_1)\ge \frac{\int_{0}^{\infty}\int_{0}^{\infty}
				\tau\,g_{1}(\tau y_{1})\,G_{2}(\tau u_1 y_{1})\,dy_{1}\,d\Pi(\tau)}
			{
				\int_{0}^{\infty}\int_{0}^{\infty} y_{1}\,\tau\,g_{1}(\tau y_{1})\,G_{2}(\tau u_1 y_{1})\,dy_{1}\,d\Pi(\tau)} = \phi_2^{*} ~~(\text{say}).
		\end{equation*}
		
		After some simplification, we get
		
		$$ \phi_2^{*}=\frac{1}{m_1+m_2}\frac{\displaystyle \int_{0}^{u_1}
			t^{\frac{m_2}{2}-1}(1+t)^{-(\frac{m_1+m_2}{2})}dt  }{E[\tau^{-1}]\displaystyle \int_{0}^{u_1} t^{\frac{m_2}{2}-1}(1+t)^{-(\frac{m_1+m_2}{2}+1)}dt \ }.$$
		
		
		\section{Estimation of  $\sigma_2^2$} \label{sec3}
		In this section, we aim to obtain  improved estimators of $\sigma_2^2$. For this purpose we consider a larger class incorporating the statistics
		$S_1$ and $S_2$.
		$$\mathcal{D}_1 = \left\lbrace \delta_{\psi_{1}}=\psi_{1}(Z_1^*)S_2 \ \ \big |  \ Z_1^* =\frac{S_1}{S_2}, \ \text{$\psi_{1}$ is positive measurable function}\right\rbrace.$$
		The risk function of $\delta_{\psi_{1}}$ under squared error loss function is
		\[
		R (\delta_{\psi_{1}},\sigma_1^{2},\sigma_2^{2})= E \left[\left(\frac{\psi_1( Z_1^*)\,S_2}{\sigma_2^2}-1\right)^2\right]
		=E\left[ E\left\lbrace\left( \frac{\psi_1(Z_1^*)\,S_2}{\sigma_2^2}-1\right)^2\big | Z_1^*\right\rbrace \right]
		.\]
		\vspace{12pt}
		To obtain a minimizer of the $R (\delta_{\psi_{1}},\sigma_1^{2},\sigma_2^{2})$, we minimize the  inner conditional risk and after simplifying, we get
		\begin{align*}
			\psi_1^Q(z_1^*)&=\frac{\frac{1}{\sigma_2^2} E\big[ S_2 \mid Z_1^* = z_1^* \big]}{\frac{1}{\sigma_2^4} E\big[ S_2^2 \mid Z_1^* = z_1^* \big]}
			=\frac{
				E\!\left[\tau^{-1}\,E\!\left(V_2 | Z_1^*=z_1^*\right)\right]
			}{E\!\left[\tau^{-2}\,E\!\left(V_2^2 | Z_1^*=z_1^*\right)\right]
			}\\
			&=\frac{1}{2}\left( 1+\frac{Z_1^*}{\eta^{-2}}\right)  \left(\frac{m_1+m_2}{2}+1 \right)^{-1} \frac{E\left[ \tau^{-1}\right] }{E\left[ \tau^{-2}\right]} \\
			&\ge  \left( 1+Z_1^*\right)  \left({m_1+m_2+2} \right)^{-1} \frac{E\left[ \tau^{-1}\right] }{E\left[ \tau^{-2}\right]} \ =\psi_{11}^Q(Z_1^*) \  (say).
		\end{align*}
		Using the convexity of the risk function we have the following theorem.  
		\begin{theorem}\label{ST12Q}
			The estimator
			$$\delta_{S1}^{Q*} = \max\left\lbrace\psi_{1}(Z_1^{*}),\psi_{11}^{Q}(Z_1^*) \right\rbrace S_2$$  has uniformly smaller risk than that of the estimator $\delta_{\psi_{1}}$ under the squared error loss function $L_1(.)$ provided $P\left(\psi^Q_{11}(Z_1^*)>\psi_1(Z_1^*)\right)>0.$
		\end{theorem}
		\noindent	Similarly, for the entropy loss function risk function is minimized at
		\begin{align*}
			\psi_1^{E}(z_1^*) &= \left(1+\frac{z_1^* }{\eta^{-2}}\right)(m_1+m_2+2)^{-1}\left(E\left[ \tau^{-1} \right] \right)^{-1} \\  &\ge \left(1+z_1^* \right)(m_1+m_2+2)^{-1}\left(E\left[\tau^{-1} \right] \right)^{-1}
			= \psi_{11}^{E}(z_1^*) \mbox{ (say). }
		\end{align*}
		In the following theorem we propose an estimator that dominates the estimator $\delta_{\psi_{1}}$ with respect to entropy loss function.
		\begin{theorem}\label{ST12E}
			The estimator
			$\delta_{S1}^{E^*} = \max\left\lbrace\psi_{1}(Z_1^*),\psi_{11}^{E}(Z_1^*) \right\rbrace S_2$ dominates $\delta_{\psi_{1}}$ with respect to the entropy loss function $L_2(.)$, provided $P\left(\psi^E_{11}(Z_1^*)>\psi_1(Z_1^*)\right)>0.$
		\end{theorem}
		\noindent Applying the Theorems \ref{ST12Q} and \ref{ST12E} we obtained the following corollary  which yields the improved estimators over BAEE. 
		\begin{corollary}
			\begin{itemize}
				\item [(i)] Under the squared error loss function the estimator 
				\begin{equation}
					\delta_{21}^Q = 
					\begin{cases}
						\frac{(1+Z_1^*)\,S_2 }{m_1 + m_2 + 2}\frac{E\left[ \tau^{-1}\right] }{E\left[ \tau^{-2}\right]}, 
						& \text{if }   Z_1^* \ge \dfrac{m_1}{m_2 + 2}, \\
						\dfrac{S_2}{m_2+ 2}\frac{E\left[ \tau^{-1}\right] }{E\left[ \tau^{-2}\right]}, 
						& \text{otherwise,}
					\end{cases}
				\end{equation}
				improves upon the  estimator $\delta_{02}^Q$.
				\item [(ii)]The estimator  \begin{equation}
					\delta_{21}^E = 
					\begin{cases}
						\dfrac{(1+Z_1^*)\,S_2}{m_1 + m_2}\left(E\left[ \tau^{-1} \right] \right)^{-1} , 
						& \text{if } Z_1^* \ge \frac{m_1}{m_2}, \\
						\frac{S_2}{m_2}\left(E\left[ \tau^{-1} \right] \right)^{-1} , 
						& \text{otherwise,}
					\end{cases}
				\end{equation}
				dominates the $\delta_{02}^E$ with respect to the entropy loss function.
			\end{itemize}
		\end{corollary}
		\noindent Next we consider a larger class of estimators utilizing the information contained in  $S_1, S_2$ and $\boldsymbol{X}_2$ as
		$$\mathcal{D}_2= \left\lbrace \delta_{\psi_{2}}= \psi_{2}(Z_1^*,Z_2^*)S_2 \big|  Z_1^*=\frac{S_1}{S_2},Z_2^* =\frac{||\boldsymbol{X}_2||^2}{S_2}, \psi_{2}  \ \text{is positive measurable function}\right\rbrace.$$
		The following theorem gives an improved estimator with respect to squared error loss function.
		\begin{theorem} \label{s2th1}
			Under the squared error loss function $L_1(.)$, risk of the estimator  $$\delta_{S2}^{Q^*} =\max\left\lbrace \psi_{21}^Q(Z_1^*,Z_2^*), \psi_{2}(Z_1^*,Z_2^*)\right\rbrace S_2 $$ is uniformly smaller than that of $\psi_{2}(Z_1^*,Z_2^*)S_2$ provided $P\left(\psi_{21}^Q(Z_1^*,Z_2^*)>\psi_1(Z_1^*,Z_2^*)\right)>0$
		\end{theorem}
		\noindent \textbf{Proof:}  The risk function of the estimator $\delta_{\psi_{2}}$ is 
		\[
		R = E \left[\left(\frac{\psi_2( Z_1^*,Z_2^*)\,S_2}{\sigma_2^2}-1\right)^2\right]
		=E\left[ E\left\lbrace\left( \frac{\psi_2(Z_1^*,Z_2^*)\,S_2}{\sigma_2^2}-1\right)^2\big | Z_1^*,Z_2^*\right\rbrace \right]
		.\]
		The minimizer of the inner condition expectation is obtained as
		\begin{equation}\label{re21}
			\psi_{2}(\eta,z_1^*,z_2^*) 	= \frac{%
				E \left[\tau^{-1}\,E\!\left(\left.\frac{\tau}{\sigma_2^2}\,S_2 \right|\ Z_1^*=z_1^*,  Z_2^*=z_2^*, K=k \right)\right]%
			}{
				E \left[\tau^{-2}\,E\!\left(\left.\frac{\tau^2}{\sigma_2^4}\,S_2^2 \right|\ Z_1^*=z_1^*, Z_2^*=z_2^*, K=k\right)\right]},
		\end{equation}
		where $K$ follows a Poisson distribution with parameter $\nu_2$.  For given $\tau>0$ and $K=k$, the probability density function of $V_2$ given  $Z_1^*$, $Z_2^*$  is proportional to
		$$ \eta^{\frac{m_1}{2}-2} v_2^{\frac{m_1+m_2+q+2k}{2}-1}e^{-\frac{v_2}{2}(1+\frac{z_1^*}{\eta^{-2}}+z_2^*)},~v_2>0.$$ Consequently we have 
		\begin{align*}\label{re23}
			\psi_{2}^Q(\eta,z_1^*,z_2^*) &= \frac{%
				E \left[\tau^{-1}\int_{0}^{\infty}v_2^{\frac{m_1+m_2+q+2k}{2}}e^{-\frac{v_2}{2}(1+\frac{z_1^*}{\eta^{-2}}+z_2^*)}dv_2\right]%
			}{%
				E \left[\tau^{-2}\int_{0}^{\infty}v_2^{\frac{m_1+m_2+q+2k}{2}+1}e^{-\frac{v_2}{2}(1+\frac{z_1^*}{\eta^{-2}}+z_2^*)}dv_2\right]%
			}\\
			&= \frac{( 1+\frac{z_1^*}{\eta^{-2}}+z_2^*)}{2}\left(\frac{m_1+m_2+q+2k}{2}+1\right)^{-1}\frac{E(\tau^{-1})}{E(\tau^{-2})}\\
			&= \frac{( 1+\frac{z_1^*}{\eta^{-2}}+z_2^*)}{m_1+m_2+q+2k+2} \frac{E(\tau^{-1})}{E(\tau^{-2})} \ge \frac{( 1+z_1^*+z_2^*)}{m_1+m_2+q+2k+2} \frac{E(\tau^{-1})}{E(\tau^{-2})}.
		\end{align*}
		\[\psi_2^Q(\eta,z_1^*,z_2^*) \ge \frac{( 1+z_1^*+z_2^*)}{m_1+m_2+q+2} \frac{E(\tau^{-1})}{E(\tau^{-2})} =\psi_{21}^Q(z_1^*,z_2^*) \ \ \text{(say)}.\]

		\noindent	In the following theorem we have given improved estimator with respect to the entropy loss function $L_2(.)$. The proof of the theorem is similar to the Theorem \ref{s2th1}, so we omit the proof. 
		\begin{theorem}\label{s2th2}
			Under the entropy loss function $L_2(.)$, risk of the estimator  $$\delta_{S2}^{E^*} =\max\left\lbrace \psi_{21}^E(Z_1^*,Z_2^*), \psi_{2}^E(Z_1^*,Z_2^*)\right\rbrace S_2$$ is nowhere larger than that of 
			$\delta_{\psi_2}^E$  provided $P\left(\psi^E_{21}(Z_1^*,Z_2^*)>\psi_2(Z_1^*,Z_2^*)\right)>0$, where $$\psi_{21}^E(Z_1^*,Z_2^*) =\frac{(1 + Z_1^* + Z_2^*)S_2}{\,m_1 + m_2 + q - 4\,}
			E[(\tau^{-1})]^{-1}.$$
		\end{theorem}
		\noindent 	Using the the Theorem \ref{s2th1} and Theorem \ref{s2th2}, we propose estimator that dominates the BAEE of $\sigma_2^2$.
		\begin{corollary}
			\begin{itemize}
				\item [(i)] The estimator 
				\begin{equation*}
					\delta_{22}^Q =
					\begin{cases}
						\frac{(1 + Z_1^* + Z_2^*)S_2}{\,m_1 + m_2 + q + 2\,}
						\;\frac{E(\tau^{-1})}{E(\tau^{-2})},
						& \text{if } 
						Z_2^* \ge \dfrac{m_1 + q}{m_2 + 2} - Z_1^*
						\ \text{ and }\ 
						Z_1^* \ge \dfrac{m_1}{m_2 + 2},
						\\
						\frac{S_2}{m_2 + 2}\;
						\frac{E(\tau^{-1})}{E(\tau^{-2})},
						& \text{otherwise},
					\end{cases}
				\end{equation*}
				improves upon the estimator $\delta_{02}^Q$ with respect to squared error loss function $L_1(.)$.
				\item [(ii)]Under the entropy loss function $L_2(.)$ the estimator  
				\[
				\delta_{22}^E=
				\begin{cases}
					\frac{(1 + Z_1^* + Z_2^*)S_2}{\,m_1 + m_2 + q - 4\,}
					\;[E(\tau^{-1})]^{-1},
					& \text{if } 
					Z_2^* \ge \dfrac{m_1 + q-4}{m_2} - Z_1^*
					\ \text{ and }\ 
					Z_1^* \ge\dfrac{m_1}{m_2},
					\\
					\frac{S_2}{m_2}\;
					[E(\tau^{-1})]^{-1},
					& \text{otherwise},
				\end{cases}
				\]
				dominates $\delta_{02}^E$.
			\end{itemize}
		\end{corollary}
		\subsection{A class of improved estimator} \label{sec3.1}
		In this section we have derived a class of estimators that dominate the BAEE. 
		\begin{theorem}
			The risk of estimator $\delta_{\psi_{1}}= {\psi_1(Z_1^*)}S_2$ is nowhere larger than $\delta_{02}^Q$ under the squared error loss function $L_1(.)$ provided the function $\psi_{1}(r)$ satisfies the following conditions:
			\begin{enumerate}
				\item[(i)]$\psi_1(r)$ is non-decreasing and $\lim_{r\rightarrow0}\psi_1(r) = {\alpha_{2Q}},$
				\item[(ii)]   $\psi_1(r) \le \frac{\int_{0}^{\infty}\int_{0}^{\infty}
					\frac{v_2}{\tau}\,g_{2}(v_{2})\,(1-G_{1}(v_{2} r))\,dv_{2}\,d\Pi(\tau)}
				{
					\int_{0}^{\infty}\int_{0}^{\infty} \frac{v_2^2}{\tau^2}\,\,g_{2}(v_2)(1-G_{1}(v_{2} r))\,dv_{2}\,d\Pi(\tau)}.$
			\end{enumerate}
		\end{theorem}
		\noindent \textbf{Proof:} Proof of this theorem is similar to the \cite{petropoulos2017estimation}.
		\begin{theorem} The risk of the estimator $\delta_{\psi_{1}}= {\psi_1(Z_1^*)}S_2$ is uniformly  smaller  than that of $\delta_{02}^{E}$ with respect to entropy  loss function $L_2(.)$ provided $\psi_{1}(r)$ satisfies the following conditions:
			\begin{enumerate}
				\item[(i)]$\psi_1(r)$ is non-decreasing and $\lim_{r\rightarrow0}\psi_1(r) = {\alpha_{2E}},$
				\item[(ii)]   $\psi_1(r) \le\frac{\int_{0}^{\infty}\int_{0}^{\infty}
					\,g_{2}(v_{2})\,(1-G_{1}(v_{2} r))\,dv_{2}\,d\Pi(\tau)}
				{\int_{0}^{\infty}\int_{0}^{\infty} \frac{v_2}{\tau}\,\,g_{2}(v_2)(1-G_{1}(v_{2} r))\,dv_{2}\,d\Pi(\tau)}.$
			\end{enumerate}
		\end{theorem}
		\noindent \textbf{Proof:} The risk difference of the estimators  $\delta_{02}^E$ and $\delta_{\psi_{1}}= {\psi_1(Z_1^*)}S_2$ is 
		%
		\begin{align*}
			RD& =E  \left[\,\frac{\psi_1(Z_1^*) S_2}{\sigma_2^2} - \ln\!\left(\frac{\psi_1(Z_1^*) S_2}{\sigma_2^2}\right) - 1 \,\right] - E \left[\,\frac{\alpha_{2E} S_2}{\sigma_2^2} - \ln\!\left(\frac{\alpha_{2E} S_2}{\sigma_2^2}\right) - 1 \,\right]\\
			&=E\left[\,\psi_1(Z_1^*) Y_2 - \ln(\psi_1(Z_1^*) Y_2) - 1 \,\right] -E  \left[ \alpha_{2E} Y_2 - \ln(\alpha_{2E} Y_2) - 1 \,\right].
		\end{align*}
		Which can be written as 
		\begin{align*}
			%
			RD	&= E\Bigg[\int_{0}^{1}
			\psi_1'\Big(\eta^{2}t\tfrac{Y_{1}}{Y_{2}}\Big)\Big(\eta^{2}Y_{1}\Big)\,
			-\frac{\psi_1'\Big(\eta^{2}t\tfrac{Y_{1}}{Y_{2}}\Big)\big(\eta^{2}Y_{1}\big)}
			{\psi_1\Big(\eta^{2}t\tfrac{Y_{1}}{Y_{2}}\Big)Y_{2}} \,dt\Bigg]\\
			&\le \int_{0}^{\infty}\!\int_{0}^{\infty}\!\int_{0}^{\infty}\!\int_{0}^{1}\!
			\psi_1'\!\Big(\frac{\eta^{2} t y_{1}}{y_{2}}\Big)
			\Bigg(1-\frac{1}{\psi_1\!\big(\frac{\eta^{2}ty_{1}}{y_2}\big)\,y_{2}}\Bigg)y_1
			\tau g_{1}(\tau y_{1})\,\tau g_{2}(\tau y_{2})
			\,dt\,dy_{2}\,dy_{1}\,d\Pi(\tau).
		\end{align*}
		Using the transformations \(r=\dfrac{t y_{1}}{y_{2}}\) and $\frac{ry_2}{t}=s$ then we obtain 
		\begin{align*}
			RD&\le \int_{0}^{\infty}\!\int_{0}^{\infty}\!\int_{0}^{\infty}\int_{ry_2}^{\infty}
			\frac{ry_2^2}{t^2}\,\psi_1'\big(\eta^{2}r\big)
			\Bigg(1-\frac{1}{\psi_1\big(\eta^{2}r\big)\,y_{2}}\Bigg)
			\tau g_{1}(\tau  \frac{ry_2}{t})\,\tau g_{2}(\tau {y_2})
			\,dt\,dr\,dy_{2}\,d\Pi(\tau).\\
			&=\int_{0}^{\infty}\int_{0}^{\infty}\int_{0}^{\infty}\int_{ry_2}^{\infty}
			x\,\psi_1'(\eta^{2}r)\Bigg(1-\frac{1}{\psi_1(\eta^{2}r)\,y_{2}}\Bigg)
			\tau g_{1}(\tau s )\tau g_{2}(\tau y_2)\,
			\frac{t}{r}\;ds\,dr \,dy_{2}\,d\Pi(\tau).\\
			&=\int_{0}^{\infty}\psi_1'(\eta^{2}r)\int_{0}^{\infty}\int_{0}^{\infty}
			\left(1-\frac{1}{\psi_1(\eta^{2}r)\,y_{2}}\right) y_{2}\,\tau\,g_{2}(\tau y_{2}) \int_{\tau ry_2}^{\infty}g_{1}(s)\;ds\,dy_{2} dr d\Pi(\tau).
		\end{align*}
		For $RD\le 0$ we must have
		\[
		\int_0^\infty\int_{0}^\infty \Big( 1 - \frac{1}{\frac{v_2}{\tau}\psi_1(\eta^{2} r)}\Big)\frac{v_2}{\tau}\, g_2( v_2)\left[ 1 - G_1(v_2r)\right] \,dv_2\,d\Pi(\tau)\le 0,
		\]
		Hence using the condition (ii) we get the theorem. 
\section{Application to the multivariate $t$ distribution}
	In this section we apply our results to multivariate $t$ distribution. If we take the  $\tau$ as $\frac{\chi^2_{\nu}}{\nu}$ in (\ref{model1}) the resulting distribution is a multivariate $t$ distribution (See \cite{petropoulos2005estimation}). Here we will give improved estimators for ordered scale parameters of two multivariate $t$ distributions. For a given $\tau>0$, let $X_{11},\dots,X_{1n_1}$ be a random sample from $N_{p}\left(\boldsymbol{\mu}_1,\frac{\sigma^2_1}{\tau}\boldsymbol{I}_p\right)$ and $X_{21},\dots,X_{2n_2}$ be a random sample from $N_{q}\left(\boldsymbol{\mu}_2,\frac{\sigma^2_1}{\tau}\boldsymbol{I}_q\right)$. Then sufficiency reduces to 
$$\overline{\boldsymbol{X}}_1\sim N_p\left(\boldsymbol{\mu}_1,\frac{\sigma_1^2}{n_1\tau} I_{p}\right) , ~~S_1 \sim\frac{ \sigma_1^2}{\tau}\chi ^2_{p(n_1-1)}~~\mbox{ and }~~\overline{\boldsymbol{X}}_2\sim N_q\left(\boldsymbol{\mu}_2,\frac{\sigma_2^2}{n_2\tau} I_{q} \right),~~S_2 \sim \frac{\sigma_2^2}{\tau}\chi ^2_{q(n_2-1)}.$$ 
Take $\boldsymbol{X}_1=\sqrt{n_1}\boldsymbol{\overline{X}}_1$ and $\boldsymbol{X}_2=\sqrt{n_2}\boldsymbol{\overline{X}}_2$ then we get the setup (\ref{model2}) with $\boldsymbol{\theta}_1=\sqrt{n_1}\boldsymbol{\mu}_1$, $\boldsymbol{\theta}_2=\sqrt{n_2}\boldsymbol{\mu}_2$, $m_1=p(n_1-1)$ and $m_2=q(n_2-1)$. Also take 
$$Z_1 =\frac{S_2}{S_1},~Z_2 =\frac{||\boldsymbol{X}_1||^2}{S_1}=\frac{||\sqrt{n_1}\boldsymbol{\overline{X}}_1||^2}{S_1},~~\mbox{and }~~Z^*_1 =\frac{S_1}{S_2},~Z^*_2 =\frac{||\boldsymbol{X}_2||^2}{S_2}=\frac{||\sqrt{n_2}\boldsymbol{\overline{X}}_2||^2}{S_2}.$$
Here $\tau$ follows Gamma$(\frac{\nu}{2}, \frac{2}{\nu})$. Employing our findings, we can get estimators for $\sigma^2_1$ as well as $\sigma^2_{2}$. Therefore the BAEE for $\sigma^2_i$,$i=1,2$ under the squared error loss  is obtained as
$d_{0i}^Q = {\frac{({\nu}-4)\nu}{4}}\frac{S_i}{2+m_i}, ~i=1,2.$ Similarly under entropy loss function BAEE for $\sigma^2_i$  can be obtained as $d^E_{0i}=  \frac{(\nu-2)\nu}{4m_i}S_i, ~i=1,2.$
Using the Theorems \ref{th1}, \ref{th2} we get the following results.
\begin{corollary}
	\begin{itemize}
		\item [(i)] Under the squared error loss function $L_1(.)$ risk of the estimator 
		\begin{equation}
			d_{11}^Q= 
			\begin{cases}
				\frac{(1+ Z_1){{({\nu}-4)\nu}}}{4({m_1+m_2}+2)  }S_1, 
				& \text{if } 0 \le Z_1 \le \dfrac{m_2}{m_1 + 2}, \\
				\frac{{{({\nu}-4)\nu}}}{4(2+m_1)}S_1, 
				& \text{otherwise.}
			\end{cases}
		\end{equation}
		is nowhere larger than that of $d_{01}^Q$.
		\item [(ii)] 	Under the entropy loss function $L_2(.)$ the estimator 
		\begin{equation}
			d_{{11}}^E = 
			\begin{cases}
				\dfrac{(1+Z_1){{({\nu}-2)\nu}}\,S_1}{4(m_1 + m_2)}, 
				& \text{if } 0 \le Z_1 \le {\frac{m_2}{m_1}} , \\
				\frac{\nu (\nu-2)}{4m_1}S_1, 
				& \text{otherwise,}
			\end{cases}
		\end{equation}
		improves upon $d_{01}^E$.
	\end{itemize}
\end{corollary}
As an application of the Theorem \ref{Th3} and \ref{Th4} we get
\begin{corollary}
	\begin{itemize}
		\item [(i)] Risk of the estimator 
		\[d_{12}^Q =
		\begin{cases}
			\frac{1 + Z_1 + Z_2}{4(m_1 + m_2 + p + 2)\,}(\nu-4)\nu S_1
			& \text{if } 
			0 \le Z_2 \le \dfrac{m_2 + p}{m_1 + 2} - Z_1
			\text{ and }
			0 \le Z_1 \le \dfrac{m_2}{m_1 + 2},
			\\
			\frac{ {{({\nu}-4)\nu}}S_1}{4(m_1 + 2)}\;
			
			& \text{otherwise},
		\end{cases}
		\]
		is nowhere larger than that of $d_{01}^Q$ under the squared error loss function $L_1(.)$. 
		\item [(ii)] Under the entropy loss function $L_2(.)$ the estimator 
		\[
		d_{12}^E =
		\begin{cases}
			\frac{(1 + Z_1 + Z_2 )S_1}{4(m_1 + m_2 + p - 4)\,}
			\; (\nu -2)\nu
			& \text{if } 
			0 \le Z_2 \le \dfrac{m_2 + p-4}{m_1} - Z_1
			\text{ and }
			0 \le Z_1 \le \dfrac{m_2}{m_1},
			\\
			\frac{\nu (\nu-2)S_1}{4m_1}
			& \text{otherwise},
		\end{cases}
		\]
		dominates $d_{01}^E$.
	\end{itemize}
\end{corollary}
Using the Theorem \ref{kth1}, \ref{kth2} we can have following result.
\begin{theorem}
	The risk of estimator $d_{\phi_{1}}= {\phi_1(Z_1)}S_1$ is nowhere larger than $d_{01}^Q$ under squared error loss function $L_1(.)$ provided $\phi_{1}(u_1)$ satisfies the following conditions:
	\begin{enumerate}
		\item[(i)]$\phi_1(u_1)$ is non-decreasing and $\lim_{u_1\rightarrow\infty}\phi_1(u_1) = \frac{1}{2+m_1}	\frac{\nu (\nu-4)}{4}.$
		\item[(ii)] $\phi(u_1) \ge \phi^*(u_1)= 	\frac{\nu (\nu-4)}{4}\frac{1}{m_1+m_2+2}\frac{E[\tau^{-1}]\displaystyle \int_{0}^{u_1}
			t^{\frac{m_2}{2}-1}(1+t)^{-(\frac{m_1+m_2}{2}+1)}dt  }{E[\tau^{-2}]\displaystyle \int_{0}^{u_1} t^{\frac{m_2}{2}-1}(1+t)^{-(\frac{m_1+m_2}{2}+2)}dt \ }.$
	\end{enumerate}
\end{theorem}
\begin{theorem}
	The estimator $d_{\phi_{1}}= {\phi_1(Z_1)}S_1$ dominates the estimator $d_{01}^{E}$ with respect to entropy loss function $L_2(.)$ provided $\phi_{1}(u_1)$ satisfies the following conditions:
	\begin{enumerate}
		\item[(i)]$\phi_1(u_1)$ is non-decreasing and $\lim_{u_1\rightarrow\infty}\phi_1(u_1) = \frac{1}{m_1}	\frac{\nu (\nu-2)}{4}.$
		\item[(ii)] $\phi(u_1) \ge 	\frac{\nu (\nu-2)}{4}
		\frac{1}{m_1+m_2}\frac{\displaystyle \int_{0}^{u_1}
			t^{\frac{m_2}{2}-1}(1+t)^{-(\frac{m_1+m_2}{2})}dt  }{E[\tau^{-1}]\displaystyle \int_{0}^{u_1} t^{\frac{m_2}{2}-1}(1+t)^{-(\frac{m_1+m_2}{2}+1)}dt \ }.$
	\end{enumerate}
\end{theorem}
Similarly, we can obtain results for $\sigma_{2}^2$ using the  section \ref{sec3}.
\begin{corollary}
	\begin{itemize}
		\item [(i)] Under the squared error loss function the estimator 
		\begin{equation}
			d_{21}^Q = 
			\begin{cases}
				\frac{(1+Z_1^*)\,S_2 }{m_1 + m_2 + 2}	\frac{\nu (\nu-4)}{4},
				& \text{if }   Z_1^* \ge \dfrac{m_1}{m_2 + 2}, \\
				\dfrac{S_2}{m_2+ 2}	\frac{\nu (\nu-4)}{4},
				& \text{otherwise,}
			\end{cases}
		\end{equation}
		improves upon the  estimator $d_{02}^Q$.
		\item [(ii)]The estimator  \begin{equation}
			d_{21}^E = 
			\begin{cases}
				\frac{2(1+Z_1^*)\,S_2}{4(m_1 + m_2)}\nu (\nu-2),
				& \text{if } Z_1^* \ge \frac{m_1}{m_2} , \\
				\frac{S_2}{m_2}	\frac{\nu (\nu-2)}{4},
				& \text{otherwise,}
			\end{cases}
		\end{equation}
		dominates the $d_{02}^E$ under the entropy loss function.
	\end{itemize}
\end{corollary}
\begin{corollary}
	\begin{itemize}
		\item [(i)] The estimator 
		\begin{equation*}
			d_{22}^Q=
			\begin{cases}
				\frac{(1 + Z_1^* + Z_2^*)S_2}{\,m_1 + m_2 + q + 2\,}
				\frac{\nu (\nu-4)}{4}
				& \text{if } 
				Z_2^* \ge \dfrac{m_1 + q}{m_2 + 2} - Z_1^*
				\ \text{ and }\ 
				Z_1^* \ge \dfrac{m_1}{m_2 + 2},
				\\
				\frac{S_2}{m_2 + 2}\;
				\frac{\nu (\nu-4)}{4},
				& \text{otherwise},
			\end{cases}
		\end{equation*}
		improves upon the estimator $d_{02}^Q$ with respect to squared error loss function $L_1(.)$.
		\item [(ii)]Under the entropy loss function $L_2(.)$ the estimator  
		\[
		d_{22}^E =
		\begin{cases}
			\frac{(1 + Z_1^* + Z_2^*)S_2}{\,m_1 + m_2 + q - 4\,}
			\frac{\nu (\nu-2)}{4}
			& \text{if } 
			Z_2^* \ge \dfrac{m_1 + q-4}{m_2} - Z_1^*
			\ \text{ and }\ 
			Z_1^* \ge \ \dfrac{m_1}{m_2},
			\\
			\frac{S_2}{m_2}\;
			\frac{\nu (\nu-2)}{4}
			& \text{otherwise},
		\end{cases}
		\]
		dominates $d_{02}^E$.
	\end{itemize}
\end{corollary}
\begin{theorem}
	The risk of estimator $d_{\psi_{1}}= {\psi_1(Z_1^*)}S_2$ is nowhere larger than $d_{02}^{Q}$ under the squared error loss function $L_1(.)$ provided the function $\psi_{1}(r)$ satisfies the following conditions:
	\begin{enumerate}
		\item[(i)]$\psi_1(r)$ is non-decreasing and $\lim_{r\rightarrow0}\psi_1(r) = \frac{1}{4(2+m_2)}{\nu}(\nu-4).$
		\item[(ii)]   $\psi_1(r) \le  \frac{\int_{0}^{\infty}\int_{0}^{\infty}
			\frac{v_2}{\tau}\,g_{2}(v_{2})\,(1-G_{1}(v_{2} r))\,dv_{2}\,d\Pi(\tau)}
		{
			\int_{0}^{\infty}\int_{0}^{\infty} \frac{v_2^2}{\tau^2}\,\,g_{2}(v_2)(1-G_{1}(v_{2} r))\,dv_{2}\,d\Pi(\tau)}.$
	\end{enumerate}
\end{theorem}
\begin{theorem}	The risk of estimator $d_{\psi_{1}}= {\psi_1(Z_1^*)}S_2$ is nowhere larger than $d_{02}^{E}$ with respect to entropy loss function $L(.)$ provided $\psi_{1}(r)$ satisfies the following conditions:
	\begin{enumerate}
		\item[(i)]$\psi_1(r)$ is non-decreasing and $\lim_{r\rightarrow0}\psi_1(r) = \frac{1}{4m_2}{\nu}({\nu}-2).$
		\item[(ii)]   $\psi_1(r) \le  \frac{\int_{0}^{\infty}\int_{0}^{\infty}
			\,g_{2}(v_{2})\,(1-G_{1}(v_{2} r))\,dv_{2}\,d\Pi(\tau)}
		{
			\int_{0}^{\infty}\int_{0}^{\infty} \frac{v_2}{\tau}\,\,g_{2}(v_2)(1-G_{1}(v_{2} r))\,dv_{2}\,d\Pi(\tau)}.$
	\end{enumerate}
	
\end{theorem}


\section{Simulation study}
\begin{figure}[H] 
	\begin{center}
		\subfigure[{$(n_1,n_2)=(5,7),(\mu_1,\mu_2)=(0,0) $}]{\includegraphics[height=4.5cm,width=6.5cm]{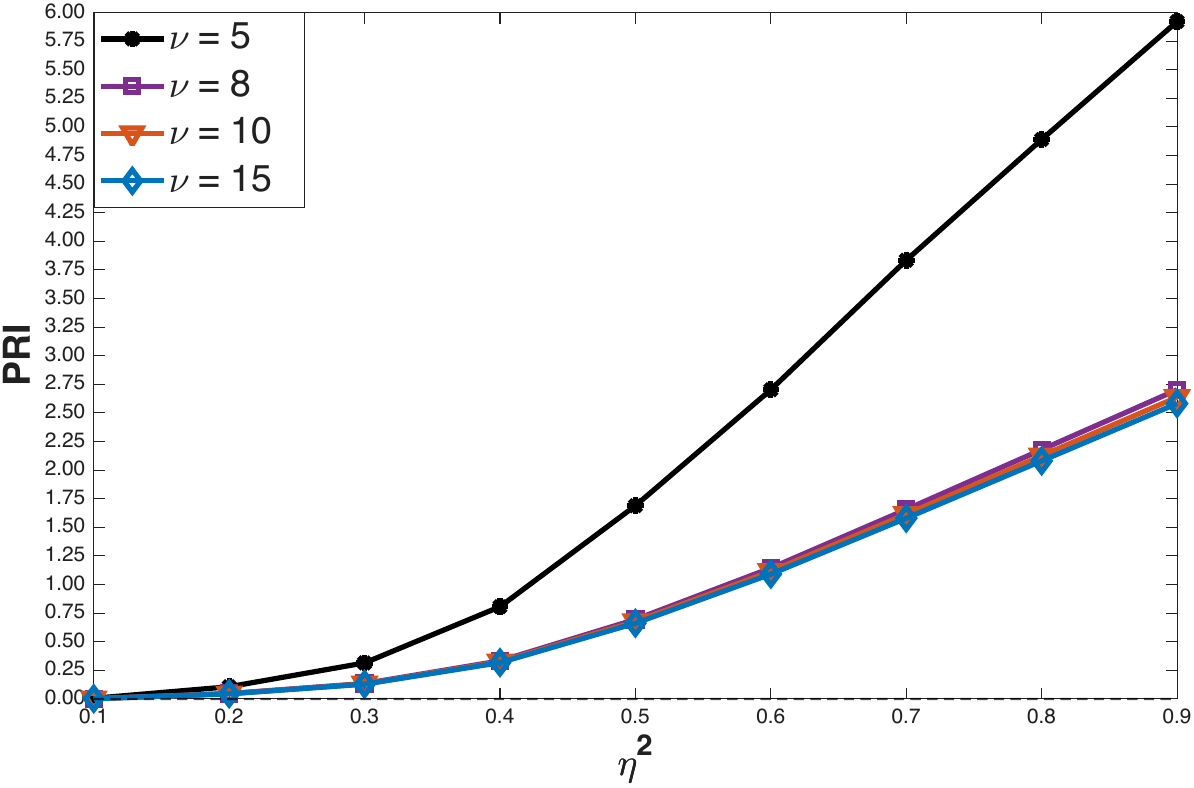}}
		\hspace{1cm} 
		\subfigure[{$(n_1,n_2)=(10,8) ,(\mu_1,\mu_2)=(0,0) $}]{\includegraphics[height=4.5cm,width=6.5cm]{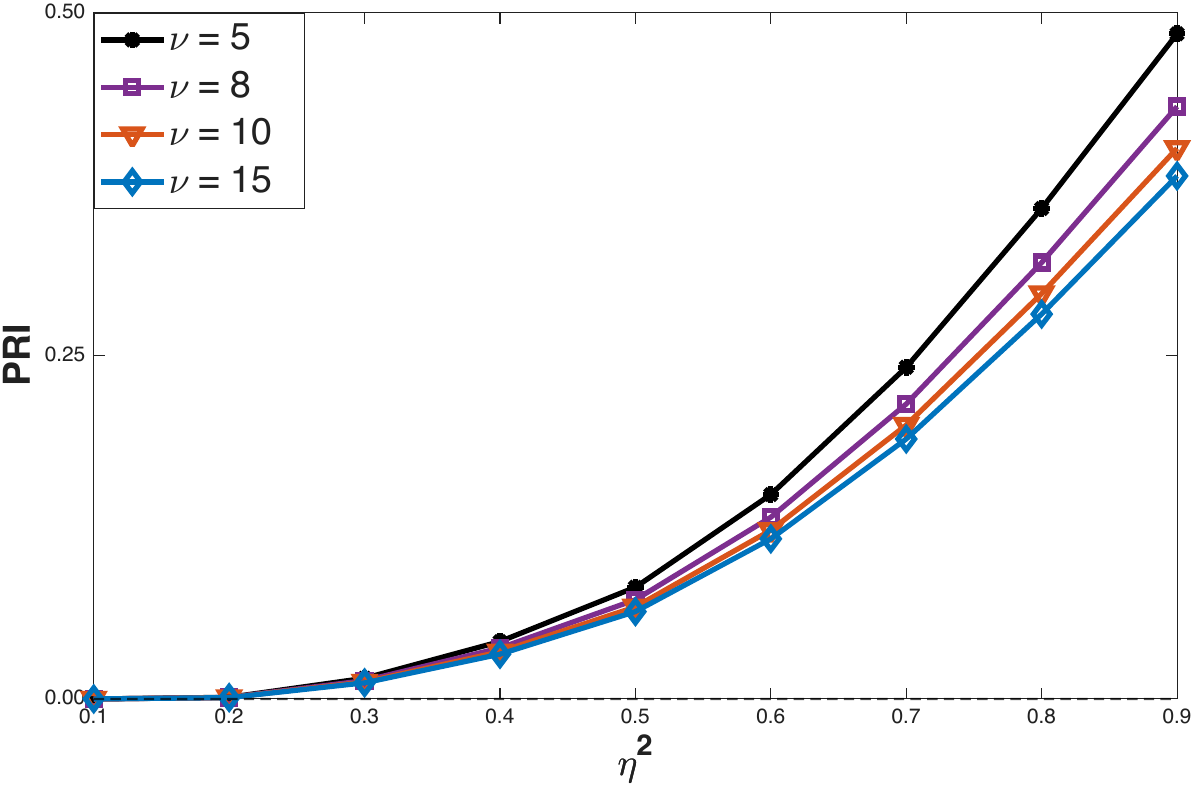}}
		\hspace{1cm}
		\subfigure[{$(n_1,n_2)=(13,9) ,(\mu_1,\mu_2)=(0,0) $}]{\includegraphics[height=4.5cm,width=6.5cm]{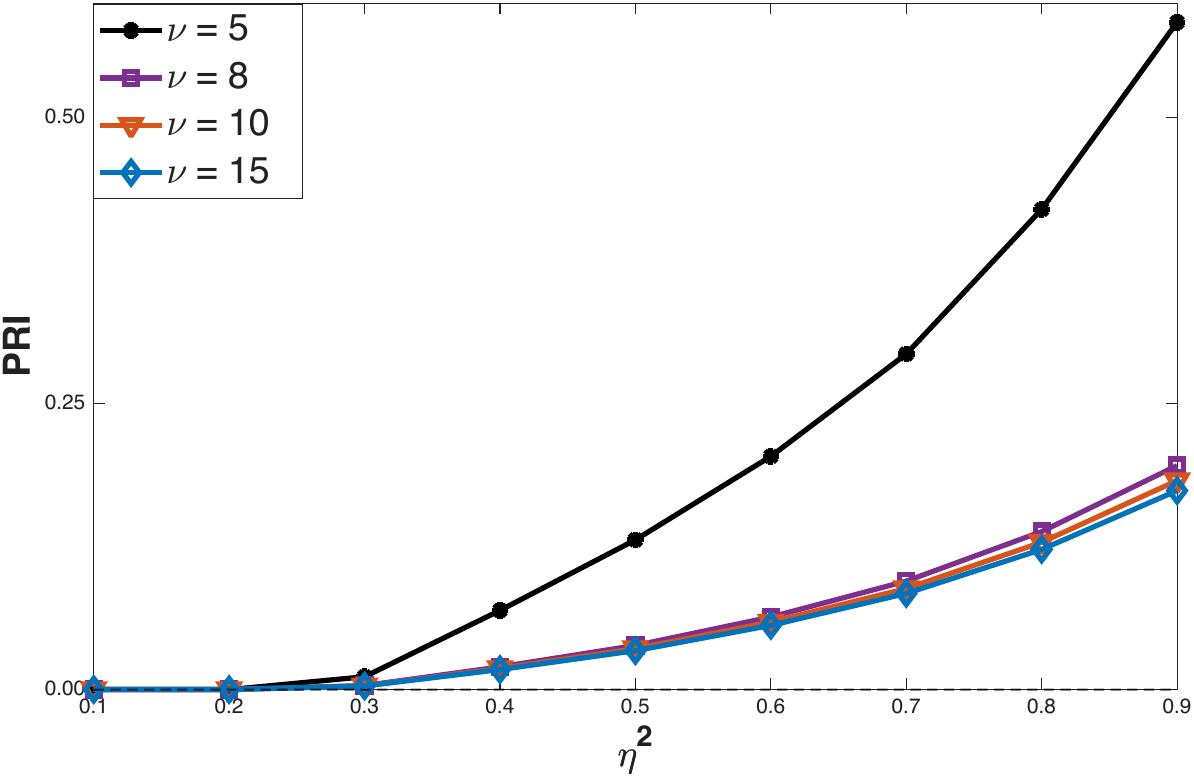}}
		\hspace{1cm}
		\subfigure[{$(n_1,n_2)=(12,15),(\mu_1,\mu_2)=(0,0)
			$}]{\includegraphics[height=4.5cm,width=6.5cm]{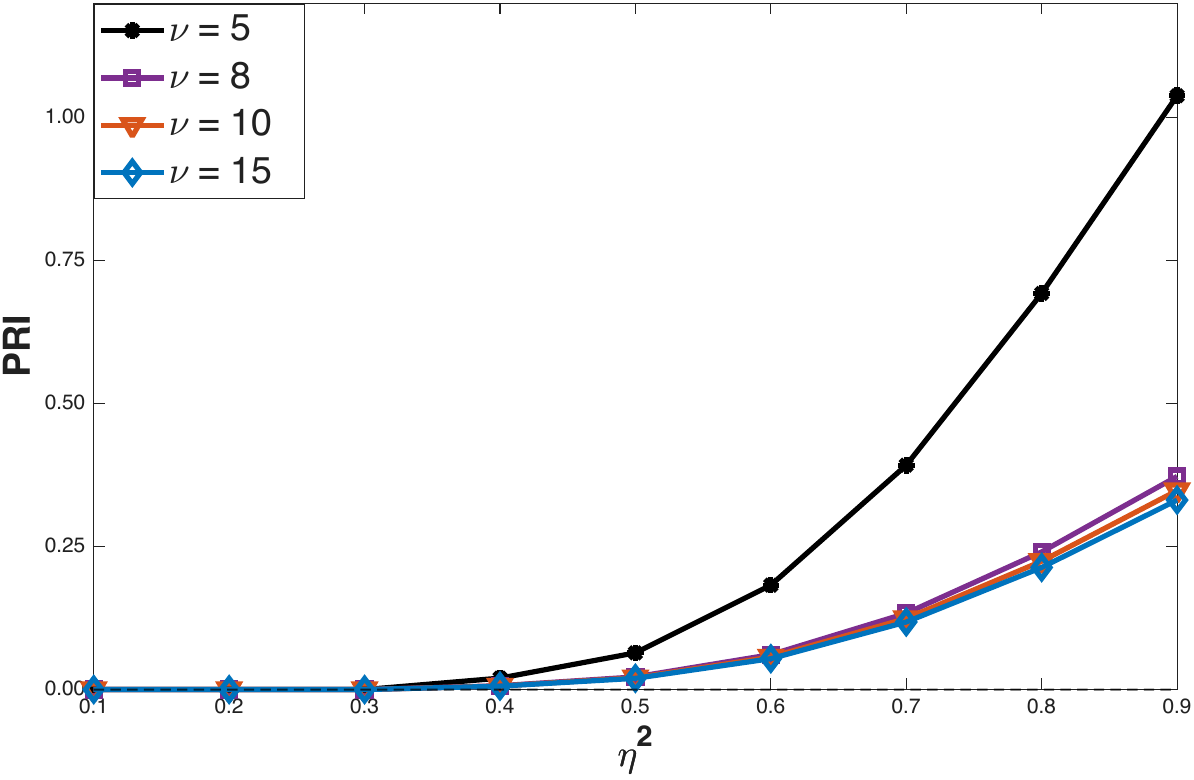}}
		\hspace{1cm}
		\subfigure[{$(n_1,n_2)=(5,7) ,(\mu_1,\mu_2)=(3,2) $}]{\includegraphics[height=4.5cm,width=6.5cm]{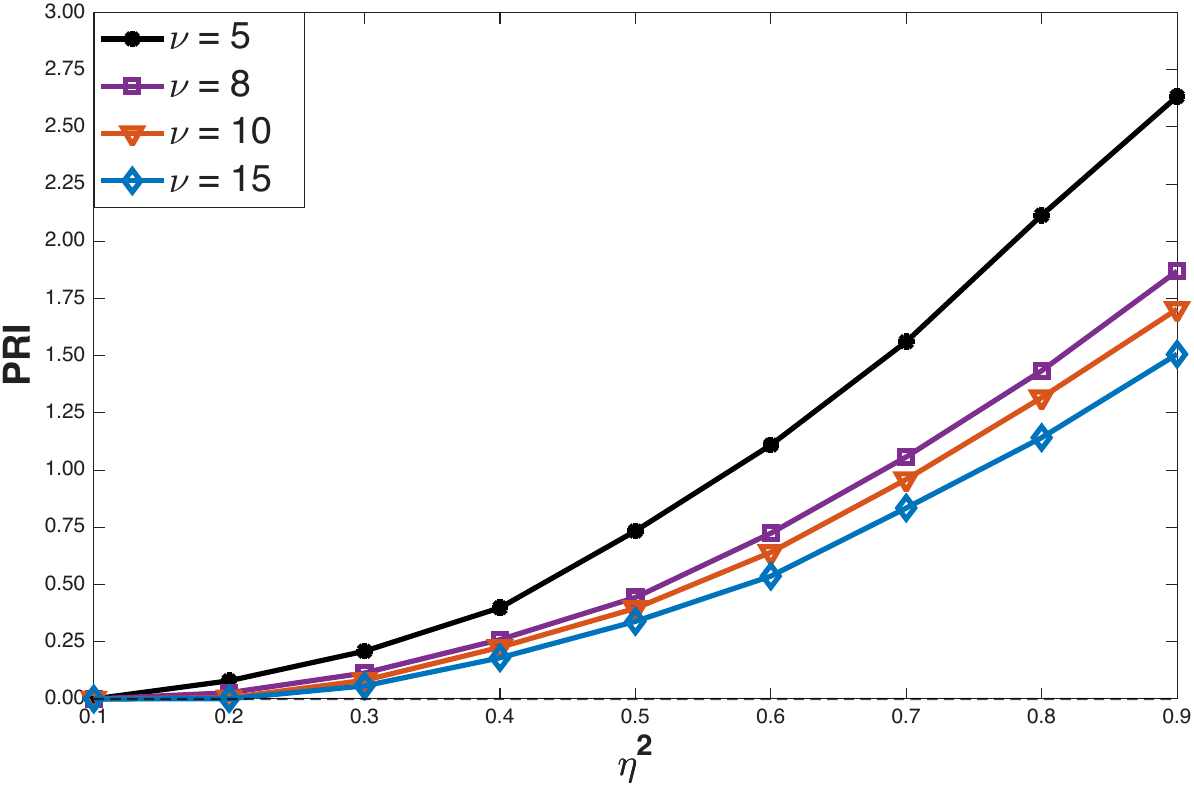}}
		\hspace{1cm} 
		\subfigure[{$(n_1,n_2)=(10,8) ,(\mu_1,\mu_2)=(3,2) $}]{\includegraphics[height=4.5cm,width=6.5cm]{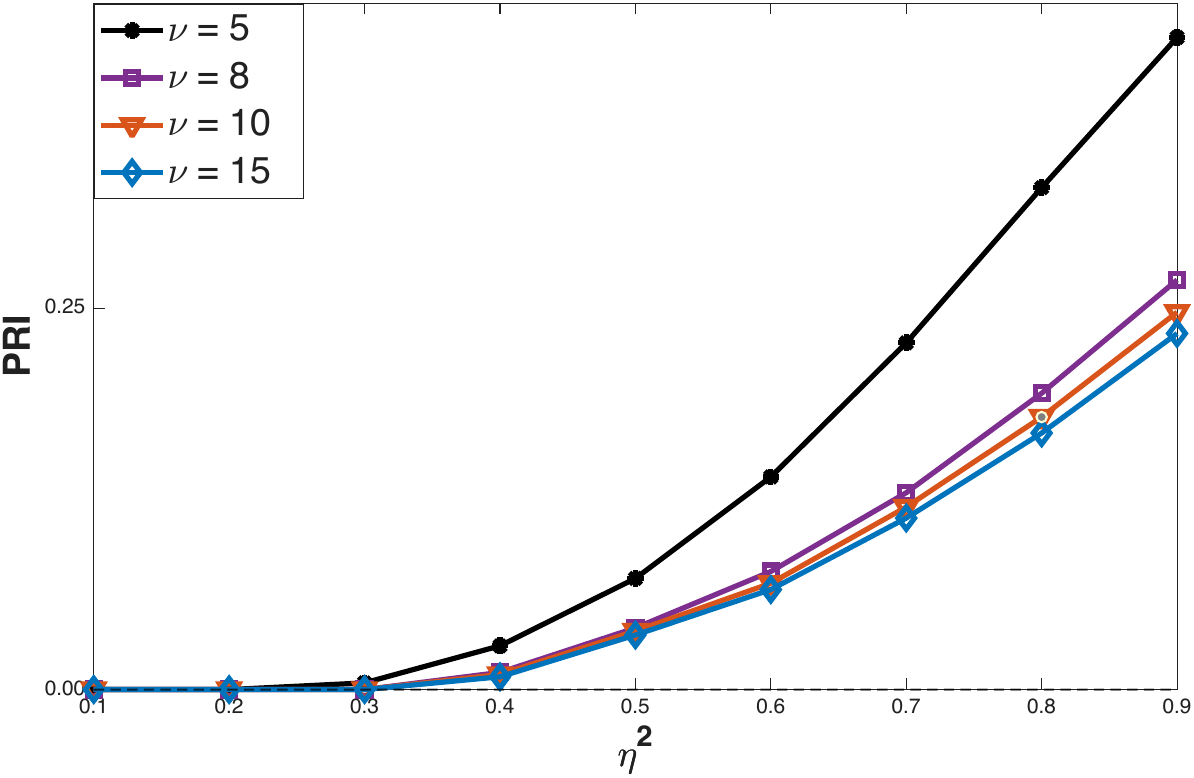}} 
		\hspace{1cm}
		\subfigure[{$(n_1,n_2)=(13,9) ,(\mu_1,\mu_2)=(3,2) $}]{\includegraphics[height=4.5cm,width=6.5cm]{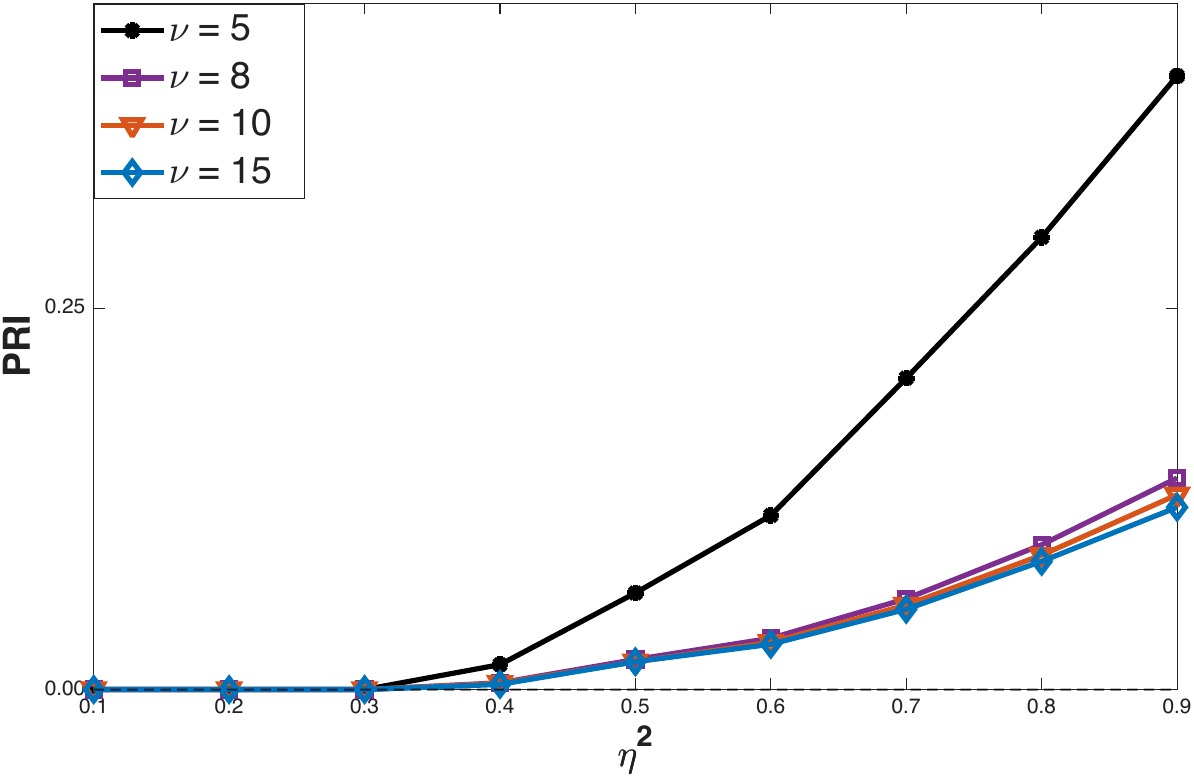}}
		\hspace{1cm} 
		\subfigure[{$(n_1,n_2)=(12,15) ,(\mu_1,\mu_2)=(3,2) $}]{\includegraphics[height=4.5cm,width=6.5cm]{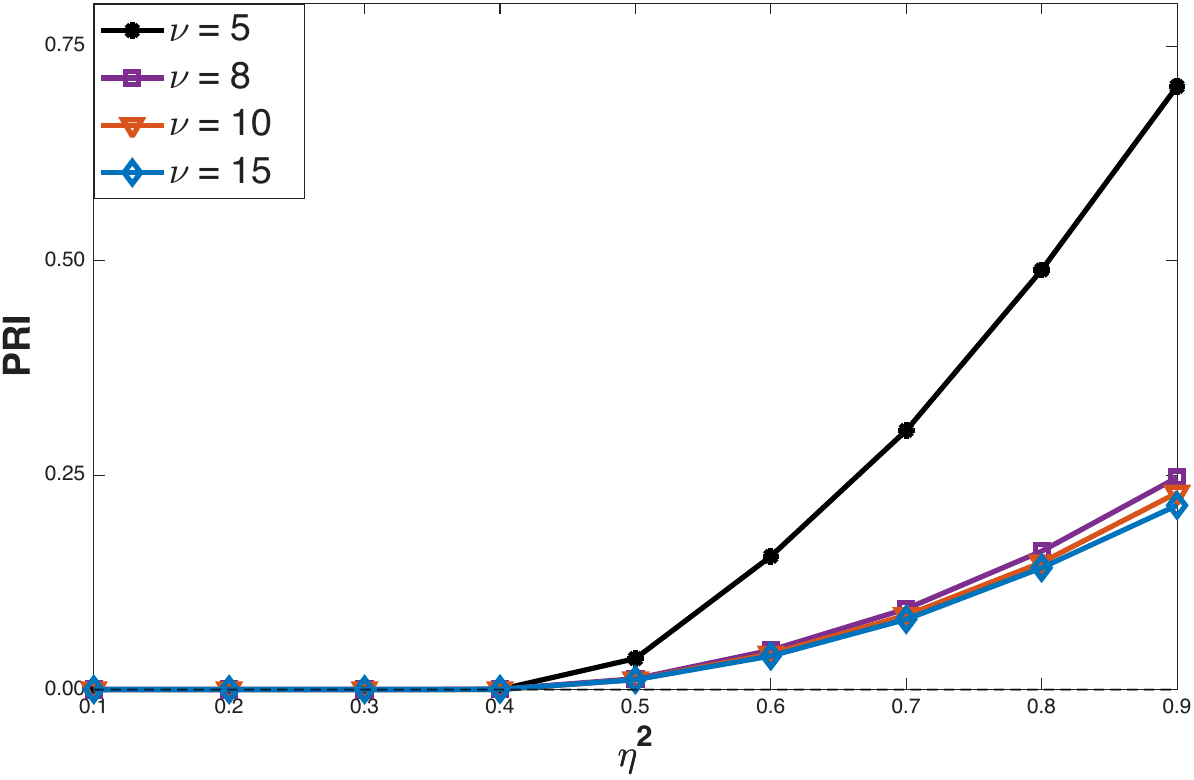}} 
		\hspace{1cm}
	\end{center}
	\caption{Percentage risk improvement of $d_{12}^Q$ under squared error loss $L_1(t)$ for $\sigma_{1}^2$}\label{figQL1}
\end{figure}

\begin{figure}[H]
	\begin{center}
		\subfigure[{$(n_1,n_2)=(5,7),(\mu_1,\mu_2)=(1,1.5) $}]{\includegraphics[height=4.5cm,width=6.5cm]{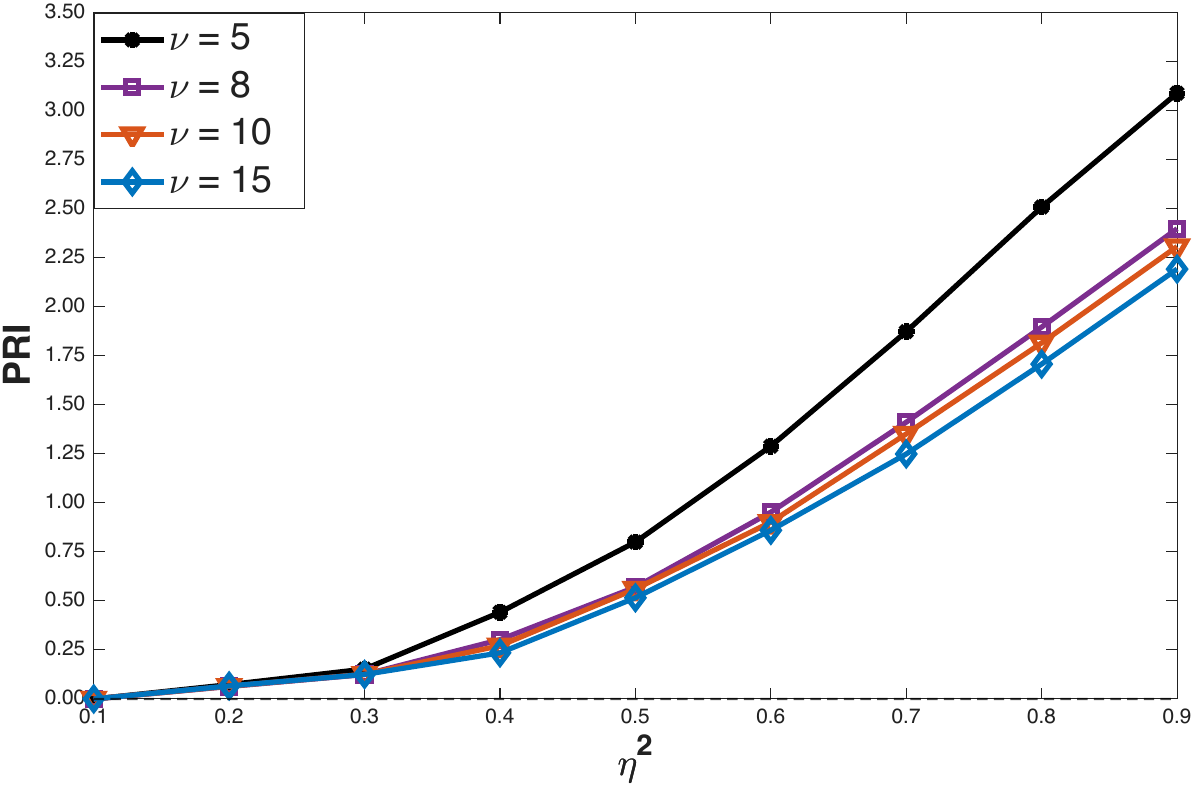}}
		\hspace{1cm} 
		\subfigure[{$(n_1,n_2)=(10,8) ,(\mu_1,\mu_2)=(1,1.5) $}]{\includegraphics[height=4.5cm,width=6.5cm]{Plots/Q_5_7_1.5.pdf}}
		\hspace{1cm}
		\subfigure[{$(n_1,n_2)=(13,9) ,(\mu_1,\mu_2)=(1,1.5) $}]{\includegraphics[height=4.5cm,width=6.5cm]{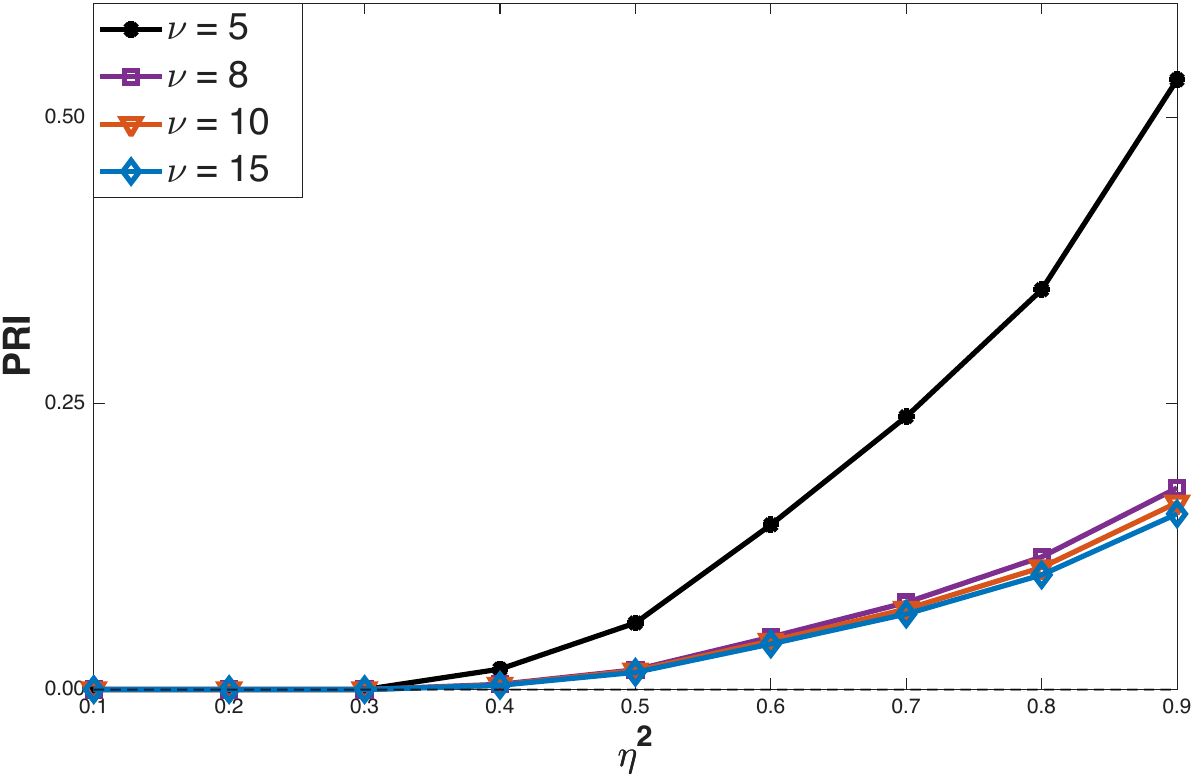}}
		\hspace{1cm}
		\subfigure[{$(n_1,n_2)=(12,15),(\mu_1,\mu_2)=(1,1.5)
			$}]{\includegraphics[height=4.5cm,width=6.5cm]{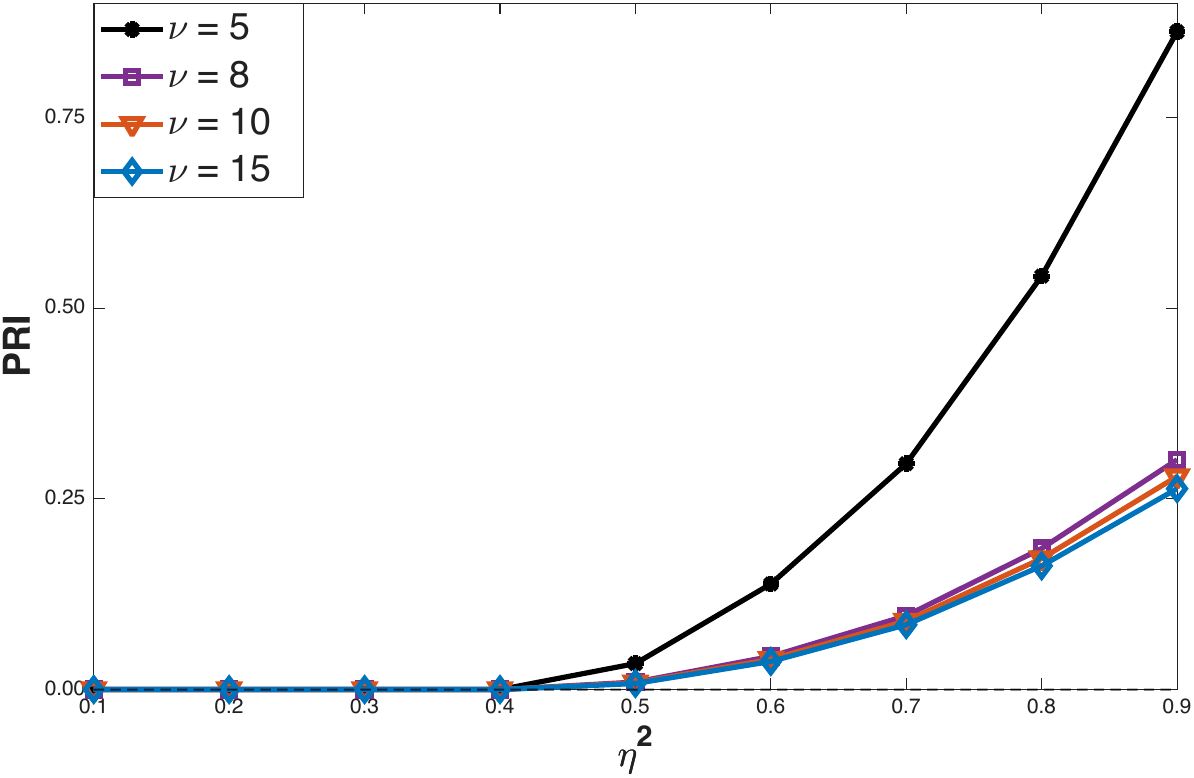}}
		\end{center}
		\caption{Percentage risk improvement of $d_{12}^Q$ under squared error loss $L_1(t)$ for $\sigma_{1}^2$}\label{figQL2}
		\end{figure}

\begin{figure}[H]
			\begin{center}
				\subfigure[{$(n_1,n_2)=(5,7),(\mu_1,\mu_2)=(0,0) $}]{\includegraphics[height=4.5cm,width=6.5cm]{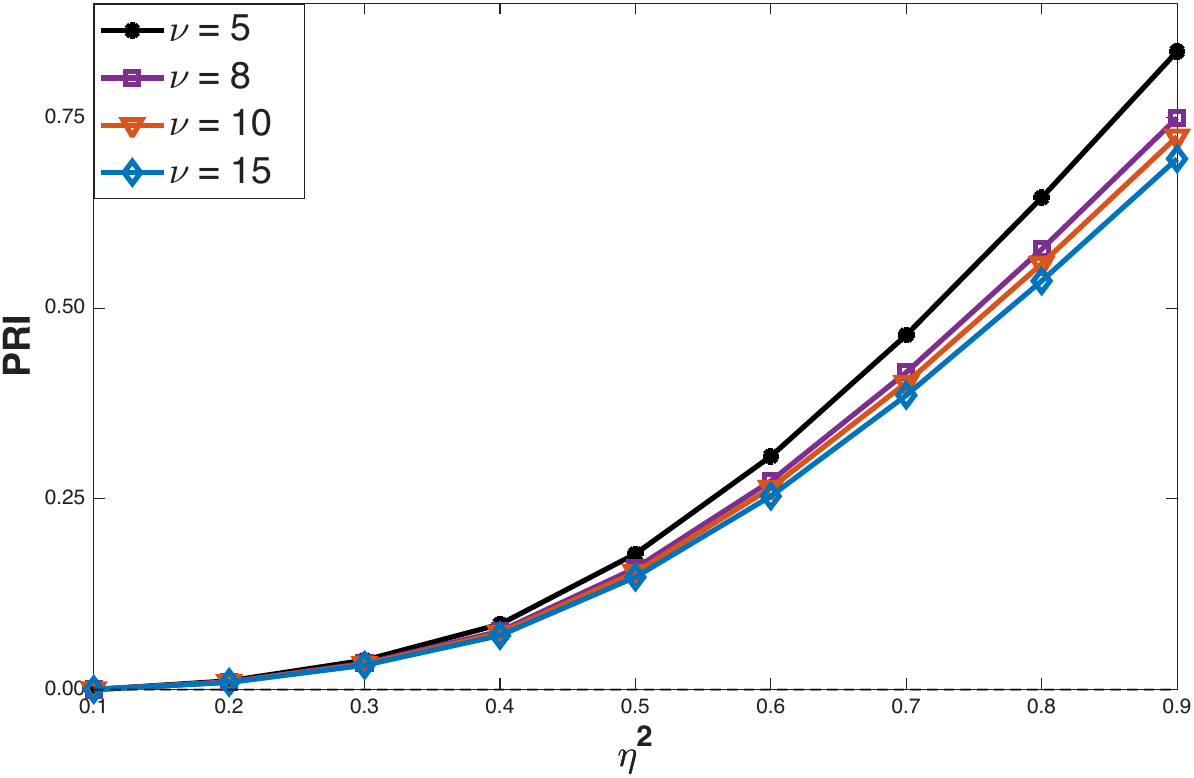}}
				\hspace{1cm} 
				\subfigure[{$(n_1,n_2)=(10,8) ,(\mu_1,\mu_2)=(0,0) $}]{\includegraphics[height=4.5cm,width=6.5cm]{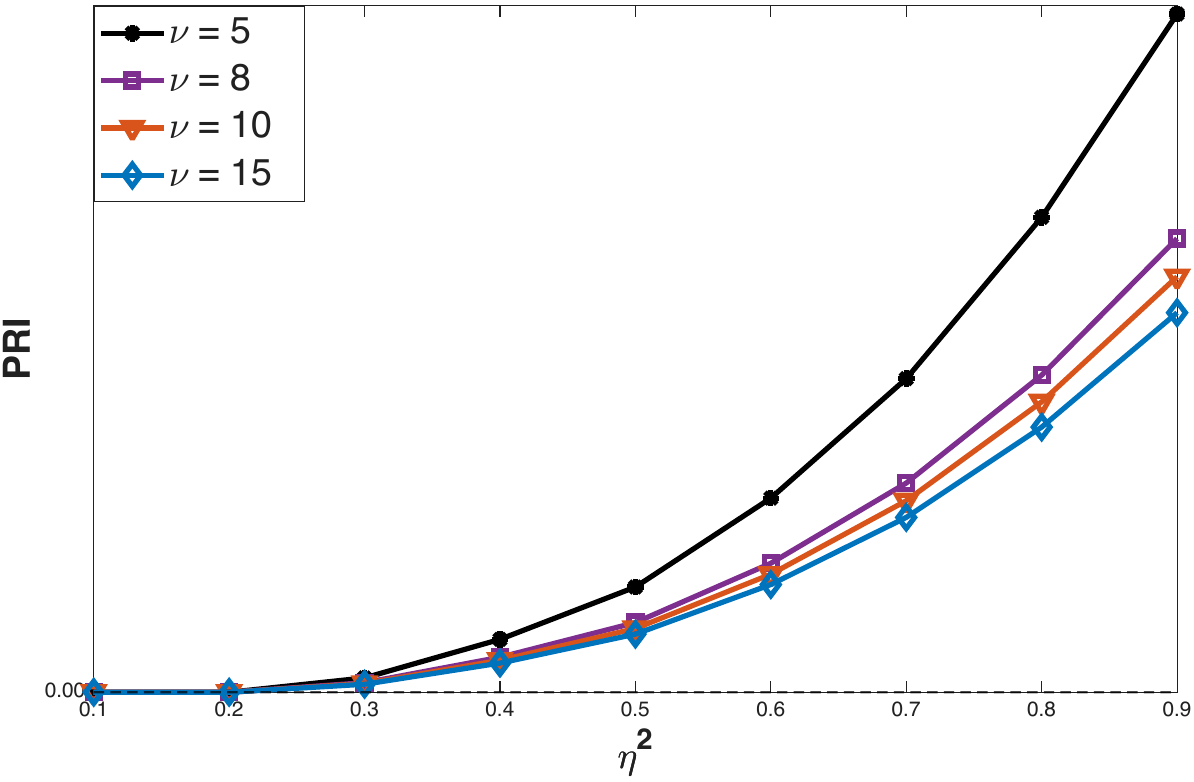}}
				\hspace{1cm}
				\subfigure[{$(n_1,n_2)=(13,9),(\mu_1,\mu_2)=(0,0) $}]{\includegraphics[height=4.5cm,width=6.5cm]{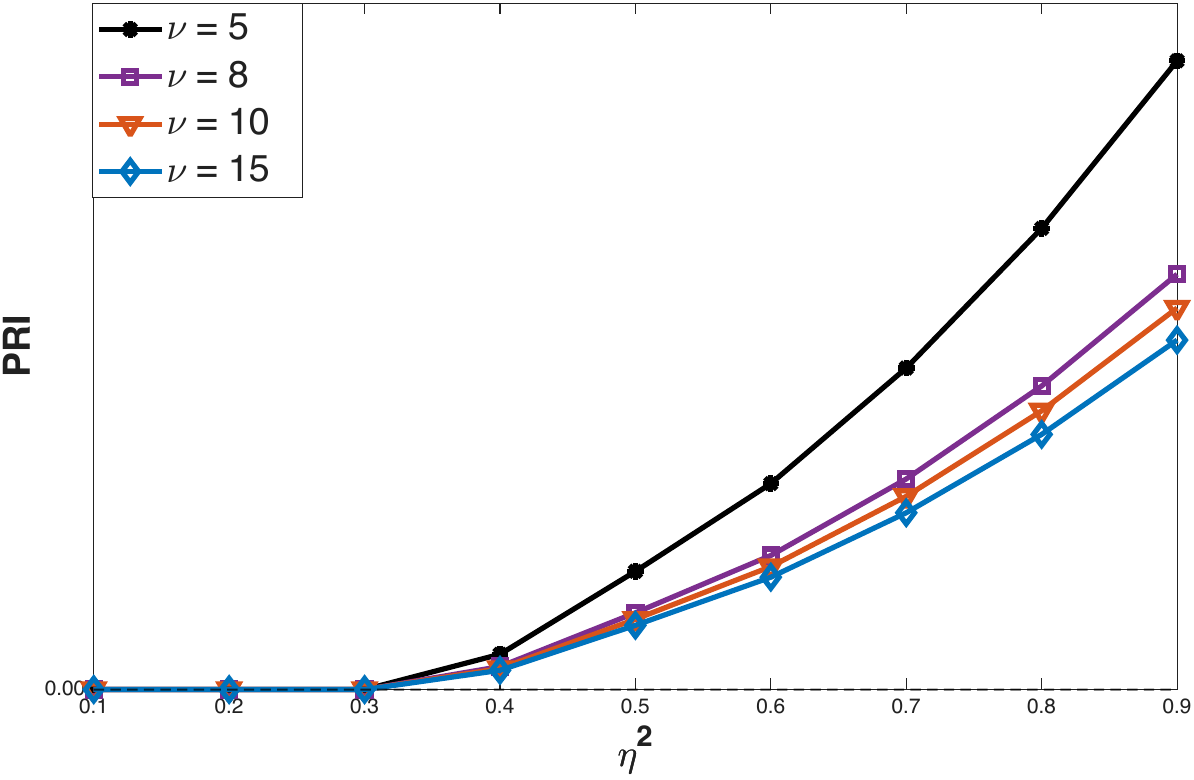}}
				\hspace{1cm} 
				\subfigure[{$(n_1,n_2)=(12,15) ,(\mu_1,\mu_2)=(0,0) $}]{\includegraphics[height=4.5cm,width=6.5cm]{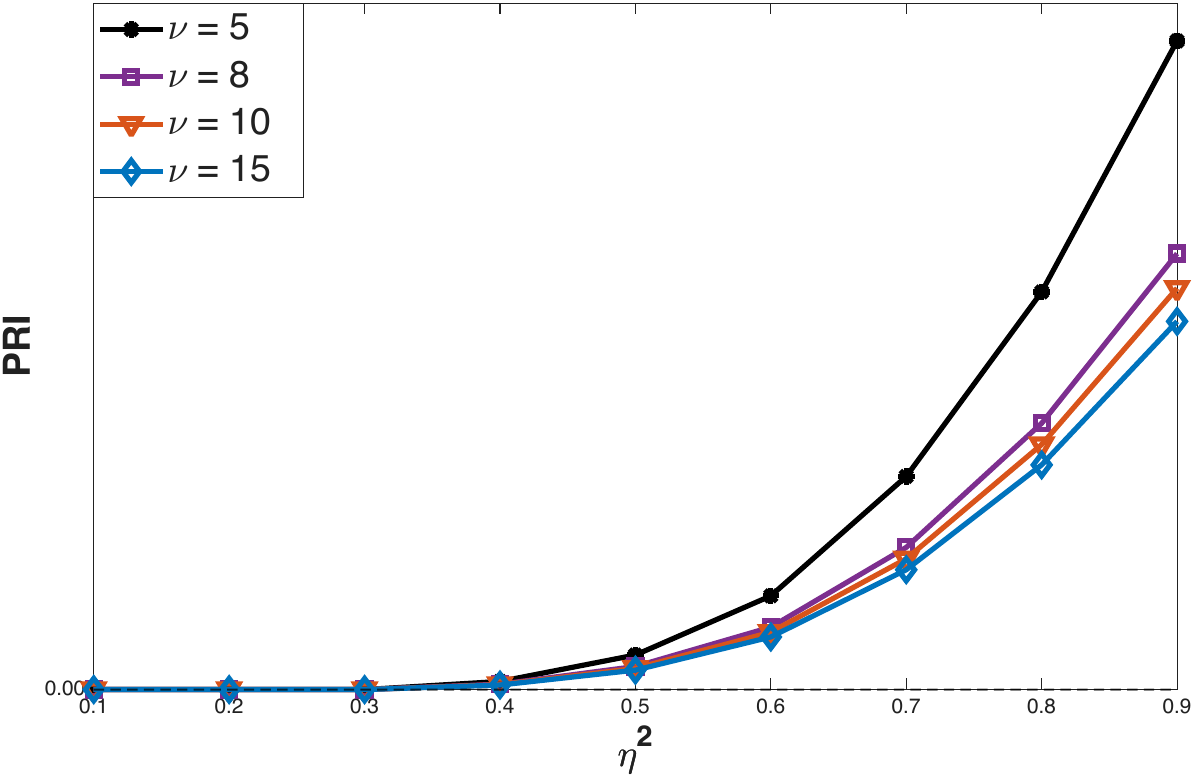}}
	\subfigure[{$(n_1,n_2)=(5,7),(\mu_1,\mu_2)=(3,2) $}]{\includegraphics[height=4.5cm,width=6.5cm]{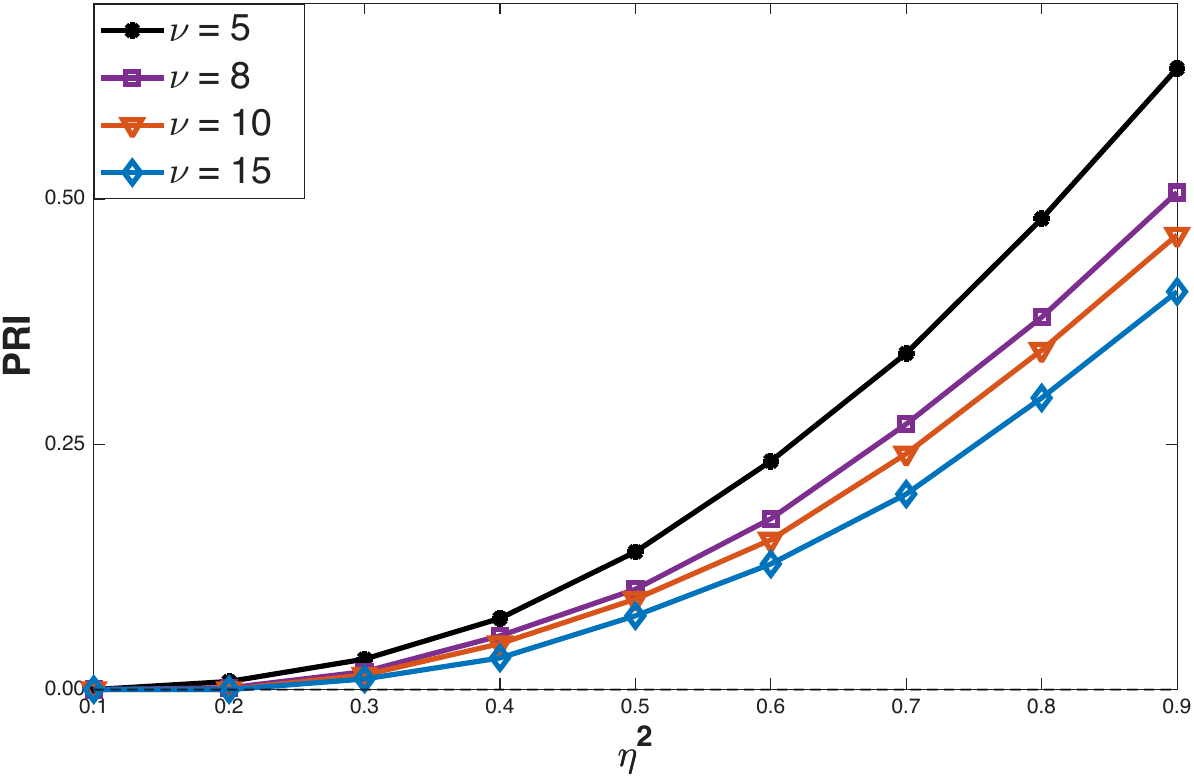}}
\hspace{1cm} 
\subfigure[{$(n_1,n_2)=(10,8) ,(\mu_1,\mu_2)=(3,2) $}]{\includegraphics[height=4.5cm,width=6.5cm]{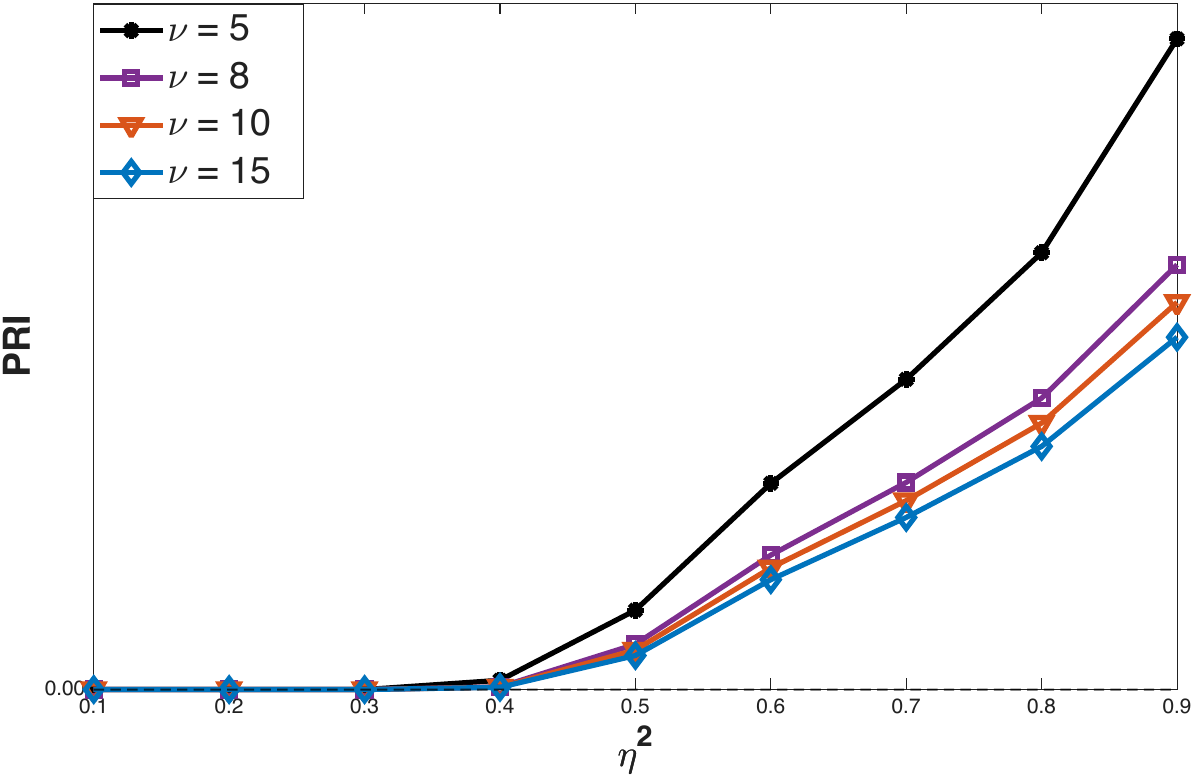}}
\hspace{1cm}
\subfigure[{$(n_1,n_2)=(13,9),(\mu_1,\mu_2)=(3,2) $}]{\includegraphics[height=4.5cm,width=6.5cm]{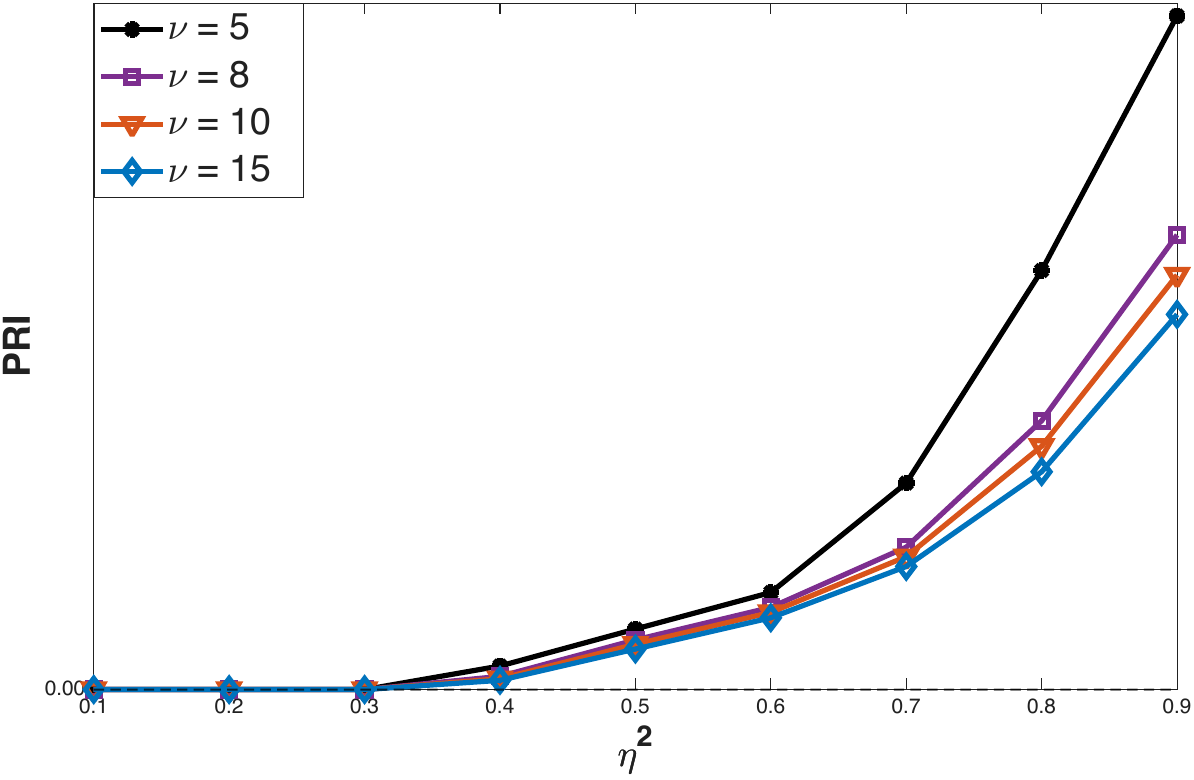}}
\hspace{1cm} 
\subfigure[{$(n_1,n_2)=(12,15) ,(\mu_1,\mu_2)=(3,2) $}]{\includegraphics[height=4.5cm,width=6.5cm]{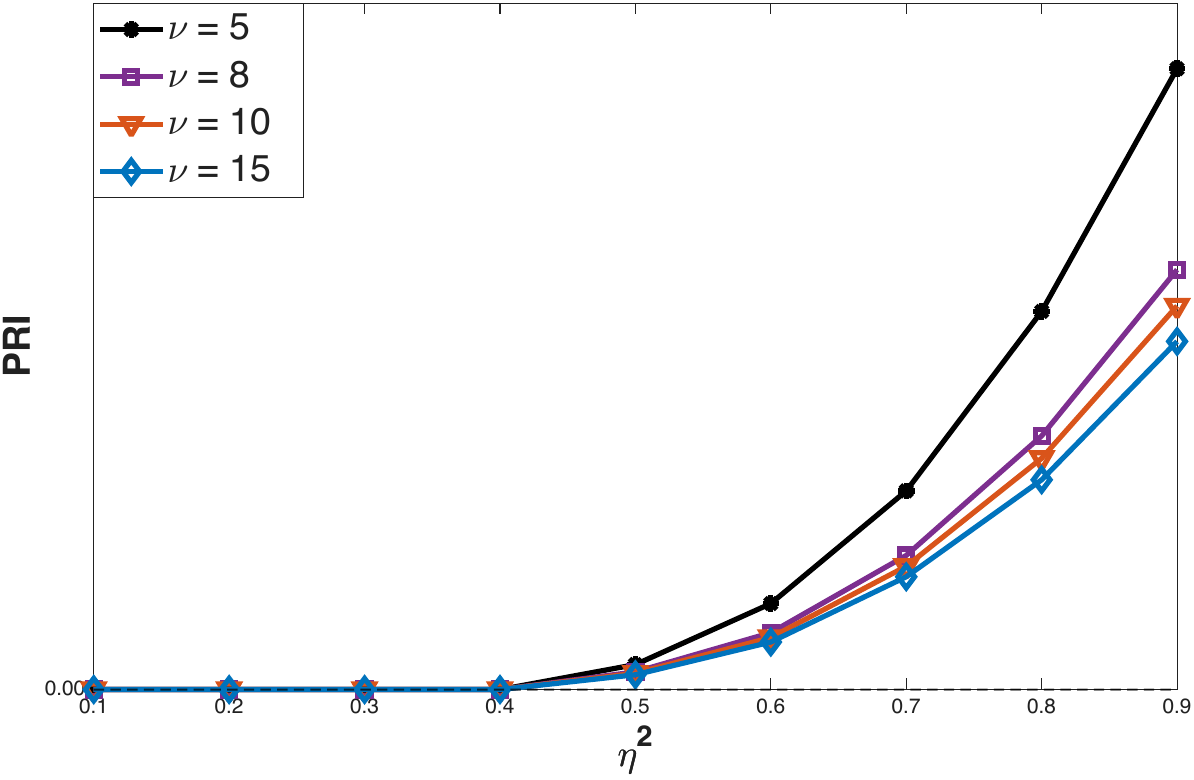}}
			\end{center}
			\caption{Percentage risk improvement of $d_{12}^E$ under entropy loss $L_2(t)$ for $\sigma_{1}^2$}\label{figEL1}
		\end{figure}

\begin{figure}[H]
\begin{center}
			\subfigure[{$(n_1,n_2)=(5,7),(\mu_1,\mu_2)=(1.5,1) $}]{\includegraphics[height=4.5cm,width=6.5cm]{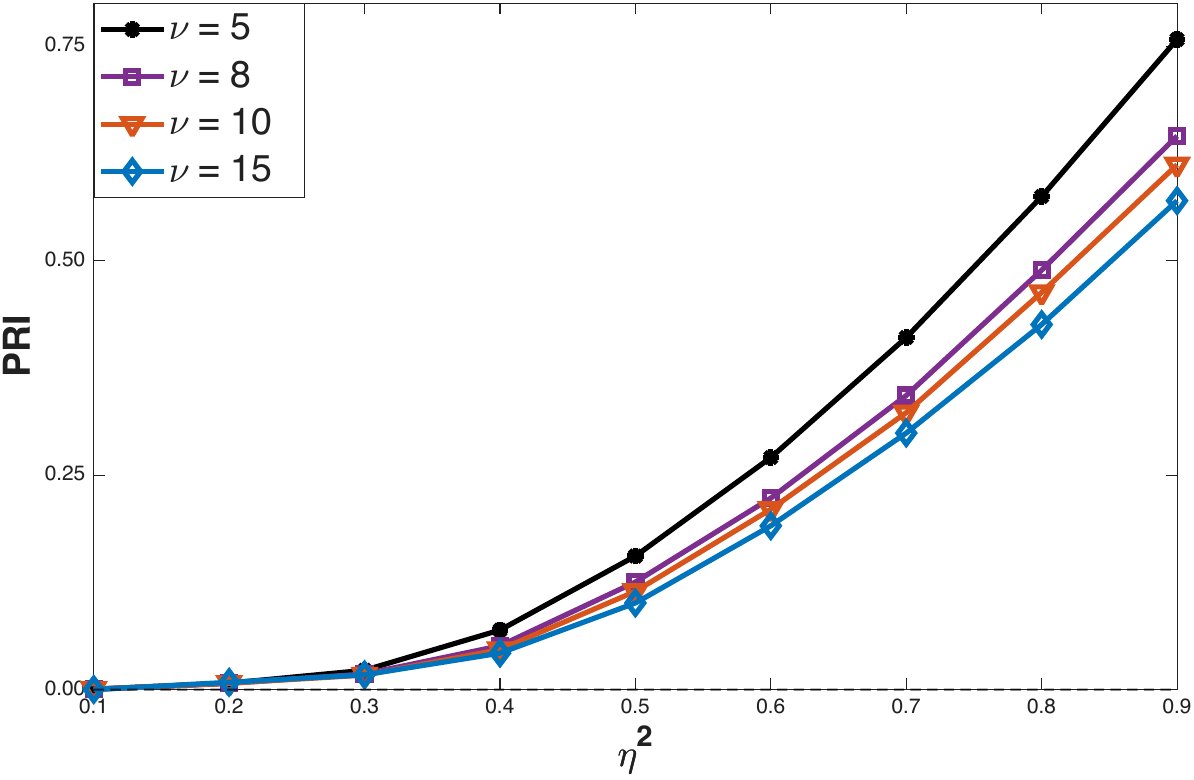}}
		\hspace{1cm} 
		\subfigure[{$(n_1,n_2)=(10,8) ,(\mu_1,\mu_2)=(1.5,1) $}]{\includegraphics[height=4.5cm,width=6.5cm]{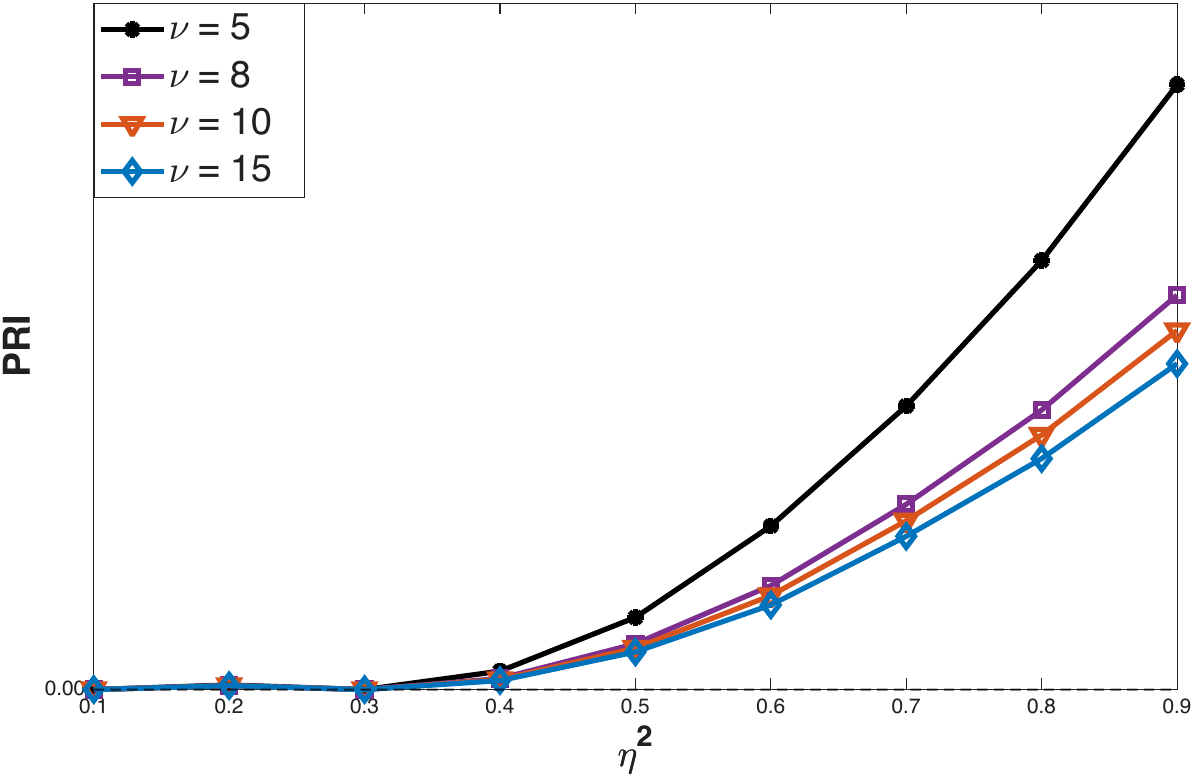}}
		\hspace{1cm}
		\subfigure[{$(n_1,n_2)=(13,9),(\mu_1,\mu_2)=(1.5,1) $}]{\includegraphics[height=4.5cm,width=6.5cm]{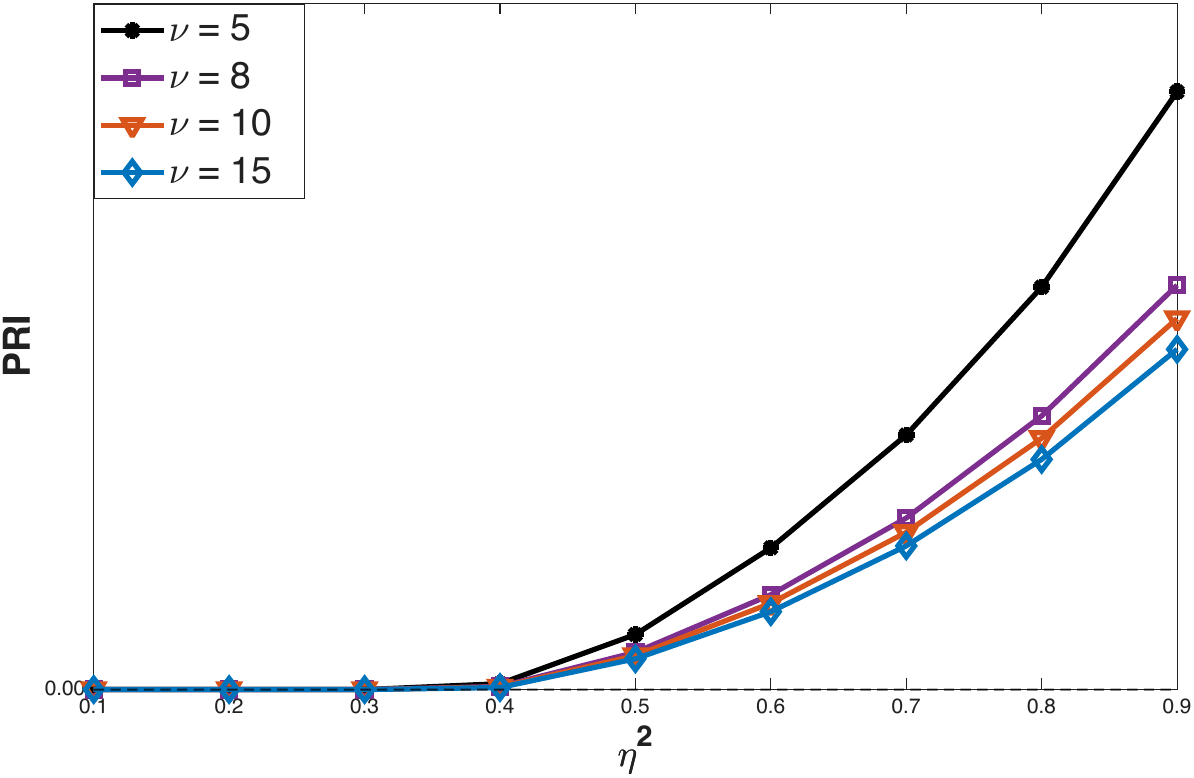}}
		\hspace{1cm} 
		\subfigure[{$(n_1,n_2)=(12,15) ,(\mu_1,\mu_2)=(1.5,1) $}]{\includegraphics[height=4.5cm,width=6.5cm]{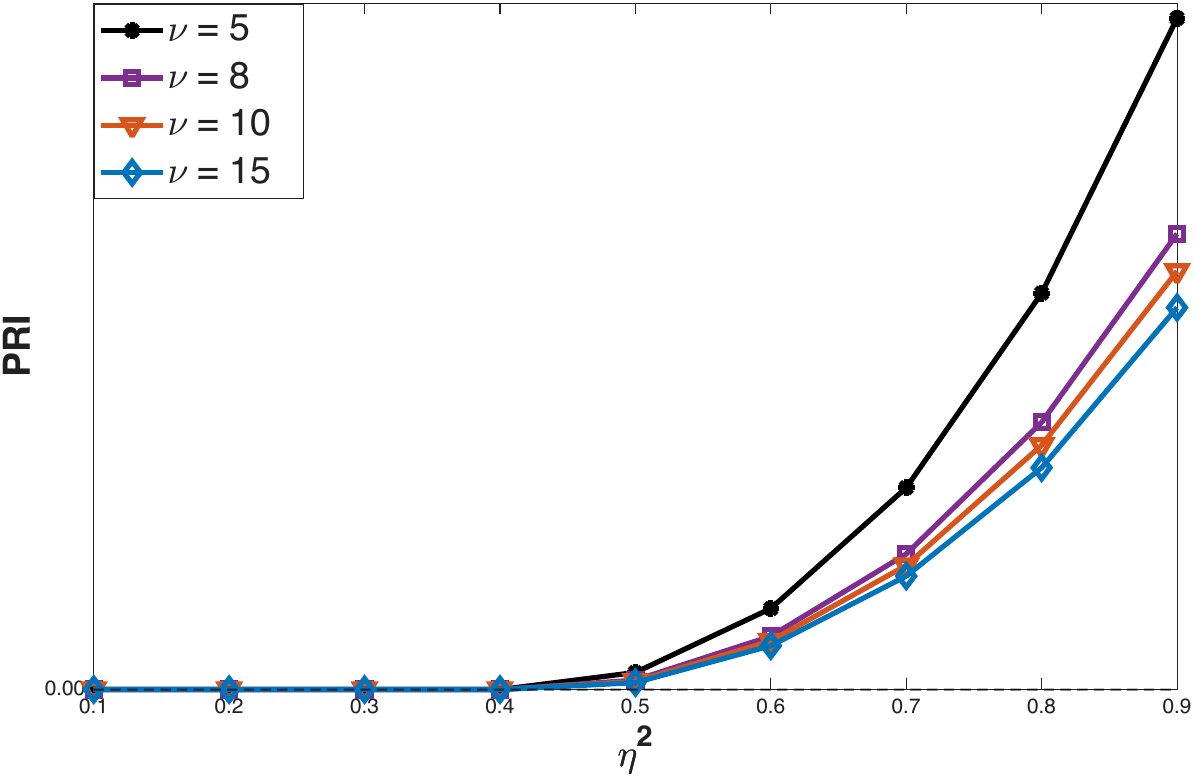}}
	\end{center}
		\caption{Percentage risk improvement of $d_{12}^E$ under entropy loss $L_2(t)$ for $\sigma_{1}^2$}\label{figEL2}
	\end{figure}

In this section, we have conducted a simulation study to asses the risk performance of the proposed estimators  for both loss functions squared error loss  and entropy loss function. Fixing $p$ and $q$ equal to 2, we have generated 20000 random samples of size $n_1$ and $n_2$ from multivariate t distribution. To compare the risk performance of the improved estimators of $\sigma_1^2$, we calculate relative  risk improvement (RRI) with respect to BAEE for different values of $\mu_1$, $\mu_2$, $\sigma_1$, $\sigma_2$ $n_1,n_2$ and $\nu$. The RRI of an estimator $\delta$  with respect to $\delta_0$ is defined as
\begin{eqnarray*}
	RRI(\delta)=\frac{Risk(\delta_0)-Risk(\delta)}{Risk(\delta_0)} \times 100.
\end{eqnarray*}
RRI of the proposed estimators $d_{11}^Q$ and $d_{BZ}^Q$  under squared error loss function is tabulated in table \ref{Tab1}. Similarly RRI $d_{11}^E$ $d_{BZ}^E$ under entropy loss presented in  table \ref{Tab2} under squared error loss function and entropy loss function respectively.
In table \ref{Tab1}.
Following insights have been drawn from table \ref{Tab1}.
\begin{enumerate}
   \item[(i)] The RRI values of both estimators $d_{11}^Q$ and $d_{BZ}^Q$ decrease monotonically as the degrees of freedom increase from $\nu = 5$ to $\nu = 15$. RRI decreases as $\frac{\sigma_{1}^2}{\sigma_{2}^2}$ increases.
   \item [(ii)] The estimator $d_{11}^Q$ yields smaller RRI values uniformly than  $d_{BZ}^Q$. This indicates $d_{BZ}^Q$ outperforms $d_{11}^Q$ under the squared error loss function. Increase in sample size smaller RRI.
\end{enumerate}
We have similar observations for entropy loss function. In the figures \ref{figQL1}, \ref{figQL2}, \ref{figEL1}, \ref{figEL2} we have demonstrated the relative risk improvement of $d_{12}^Q$ for different values of $\nu$. We have observed that smaller value of $\nu$ gives better improvement over BAEE for  $d_{12}^Q$ and  $d_{12}^E$.
\newpage
\begin{table}[h!]
	\begin{center}
		\caption{Percentage risk improvement of $d_{11}^Q$ and $d_{BZ}^Q$ with respect to squared error loss function.}\label{Tab1}
			\begin{tabular}{|c|c|cc|cc|cc|cc|}
				\hline 
				\multirow{2}{*}{$(n_1,n_2)$} & \multirow{2}{*}{$\frac{\sigma_1^2}{\sigma_2^2}$} & \multicolumn{2}{c|}{$\nu=5$}                          & \multicolumn{2}{c|}{$\nu=8$}                 & \multicolumn{2}{c|}{$\nu=10$}                                  & \multicolumn{2}{c|}{$\nu=15$}                                  \\ \cline{3-10} 
				&                                                  & \multicolumn{1}{c|}{$d_{11}^Q$} & \textbf{$d_{BZ}^Q$} & \multicolumn{1}{c|}{$d_{11}^Q$} & $d_{BZ}^Q$ & \multicolumn{1}{c|}{$d_{11}^Q$} & $d_{BZ}^Q$ & \multicolumn{1}{c|}{$d_{11}^Q$} & $d_{BZ}^Q$\\ \hline
				\multirow{6}{*}{$(5,7)$}     & 0.1                                              & \multicolumn{1}{c|}{0.02}       & 1.56                & \multicolumn{1}{c|}{0.02}       & 1.23       & \multicolumn{1}{c|}{0.02}       & 1.20                         & \multicolumn{1}{c|}{0.02}       & 1.17                         \\ \cline{2-10} 
				& 0.2                                              & \multicolumn{1}{c|}{0.40}       & 7.30                & \multicolumn{1}{c|}{0.31}       & 5.84       & \multicolumn{1}{c|}{0.30}       & 5.70                         & \multicolumn{1}{c|}{0.30}       & 5.59                         \\ \cline{2-10} 
				& 0.4                                              & \multicolumn{1}{c|}{3.59}       & 22.13               & \multicolumn{1}{c|}{2.84}       & 18.09      & \multicolumn{1}{c|}{2.77}       & 17.69                        & \multicolumn{1}{c|}{2.71}       & 17.35                        \\ \cline{2-10} 
				& 0.6                                              & \multicolumn{1}{c|}{10.91}      & 34.78               & \multicolumn{1}{c|}{8.73}       & 28.84      & \multicolumn{1}{c|}{8.52}       & 28.24                        & \multicolumn{1}{c|}{8.34}       & 27.74                        \\ \cline{2-10} 
				& 0.7                                              & \multicolumn{1}{c|}{15.39}      & 39.91               & \multicolumn{1}{c|}{12.39}      & 33.29      & \multicolumn{1}{c|}{12.10}      & 32.62                        & \multicolumn{1}{c|}{11.85}      & 32.06                        \\ \cline{2-10} 
				& 0.8                                              & \multicolumn{1}{c|}{20.03}      & 44.35               & \multicolumn{1}{c|}{16.22}      & 37.18      & \multicolumn{1}{c|}{15.85}      & 36.46                        & \multicolumn{1}{c|}{15.53}      & 35.84                        \\ \hline
				\multirow{6}{*}{$(10,8)$}    & 0.1                                              & \multicolumn{1}{c|}{0.00}       & 0.17                & \multicolumn{1}{c|}{0.00}       & 0.15       & \multicolumn{1}{c|}{0.00}       & 0.14                         & \multicolumn{1}{c|}{0.00}       & 0.13                         \\ \cline{2-10} 
				& 0.2                                              & \multicolumn{1}{c|}{0.07}       & 1.79                & \multicolumn{1}{c|}{0.06}       & 1.58       & \multicolumn{1}{c|}{0.05}       & 1.47                         & \multicolumn{1}{c|}{0.05}       & 1.39                         \\ \cline{2-10} 
				& 0.4                                              & \multicolumn{1}{c|}{1.24}       & 9.93                & \multicolumn{1}{c|}{1.10}       & 8.89       & \multicolumn{1}{c|}{1.01}       & 8.28                         & \multicolumn{1}{c|}{0.96}       & 7.88                         \\ \cline{2-10} 
				& 0.6                                              & \multicolumn{1}{c|}{5.14}       & 19.81               & \multicolumn{1}{c|}{4.56}       & 17.84      & \multicolumn{1}{c|}{4.24}       & 16.69                        & \multicolumn{1}{c|}{4.02}       & 15.92                        \\ \cline{2-10} 
				& 0.7                                              & \multicolumn{1}{c|}{8.07}       & 24.35               & \multicolumn{1}{c|}{7.19}       & 21.99      & \multicolumn{1}{c|}{6.68}       & 20.61                        & \multicolumn{1}{c|}{6.35}       & 19.68                        \\ \cline{2-10} 
				& 0.8                                              & \multicolumn{1}{c|}{11.49}      & 28.48               & \multicolumn{1}{c|}{10.26}      & 25.78      & \multicolumn{1}{c|}{9.56}       & 24.19                        & \multicolumn{1}{c|}{9.09}       & 23.13                        \\ \hline
				\multirow{6}{*}{$(12,15)$}   & 0.1                                              & \multicolumn{1}{c|}{0.00}       & 0.03                & \multicolumn{1}{c|}{0.00}       & 0.01       & \multicolumn{1}{c|}{0.00}       & 0.01                         & \multicolumn{1}{c|}{0.00}       & 0.01                         \\ \cline{2-10} 
				& 0.2                                              & \multicolumn{1}{c|}{0.00}       & 1.13                & \multicolumn{1}{c|}{0.00}       & 0.36       & \multicolumn{1}{c|}{0.00}       & 0.33                         & \multicolumn{1}{c|}{0.00}       & 0.32                         \\ \cline{2-10} 
				& 0.4                                              & \multicolumn{1}{c|}{0.89}       & 14.32               & \multicolumn{1}{c|}{0.27}       & 5.24       & \multicolumn{1}{c|}{0.25}       & 4.89                         & \multicolumn{1}{c|}{0.23}       & 4.66                         \\ \cline{2-10} 
				& 0.6                                              & \multicolumn{1}{c|}{7.43}       & 35.65               & \multicolumn{1}{c|}{2.48}       & 14.71      & \multicolumn{1}{c|}{2.31}       & 13.81                        & \multicolumn{1}{c|}{2.19}       & 13.20                        \\ \cline{2-10} 
				& 0.7                                              & \multicolumn{1}{c|}{14.36}      & 45.51               & \multicolumn{1}{c|}{5.04}       & 19.84      & \multicolumn{1}{c|}{4.70}       & 18.66                        & \multicolumn{1}{c|}{4.48}       & 17.85                        \\ \cline{2-10} 
				& 0.8                                              & \multicolumn{1}{c|}{23.00}      & 53.95               & \multicolumn{1}{c|}{8.50}       & 24.77      & \multicolumn{1}{c|}{7.94}       & 23.34                        & \multicolumn{1}{c|}{7.57}       & 22.36                        \\ \hline
				\multirow{6}{*}{$(13,18)$}   & 0.1                                              & \multicolumn{1}{c|}{0.00}       & 0.01                & \multicolumn{1}{c|}{0.00}       & 0.00       & \multicolumn{1}{c|}{0.00}       & 0.00                         & \multicolumn{1}{c|}{0.00}       & 0.00                         \\ \cline{2-10} 
				& 0.2                                              & \multicolumn{1}{c|}{0.00}       & 0.69                & \multicolumn{1}{c|}{0.00}       & 0.23       & \multicolumn{1}{c|}{0.00}       & 0.20                         & \multicolumn{1}{c|}{0.00}       & 0.18                         \\ \cline{2-10} 
				& 0.4                                              & \multicolumn{1}{c|}{0.63}       & 12.06               & \multicolumn{1}{c|}{0.19}       & 4.58       & \multicolumn{1}{c|}{0.17}       & 4.08                         & \multicolumn{1}{c|}{0.16}       & 3.79                         \\ \cline{2-10} 
				& 0.6                                              & \multicolumn{1}{c|}{6.26}       & 33.83               & \multicolumn{1}{c|}{2.17}       & 14.47      & \multicolumn{1}{c|}{1.92}       & 13.01                        & \multicolumn{1}{c|}{1.78}       & 12.13                        \\ \cline{2-10} 
				& 0.7                                              & \multicolumn{1}{c|}{13.02}      & 44.51               & \multicolumn{1}{c|}{4.77}       & 20.12      & \multicolumn{1}{c|}{4.24}       & 18.15                        & \multicolumn{1}{c|}{3.93}       & 16.96                        \\ \cline{2-10} 
				& 0.8                                              & \multicolumn{1}{c|}{21.95}      & 53.83               & \multicolumn{1}{c|}{8.47}       & 25.65      & \multicolumn{1}{c|}{7.56}       & 23.22                        & \multicolumn{1}{c|}{7.02}       & 21.74                        \\ \hline
				\multirow{6}{*}{$(20,17)$}   & 0.1                                              & \multicolumn{1}{c|}{0.00}       & 0.00                & \multicolumn{1}{c|}{0}          & 0.00       & \multicolumn{1}{c|}{0.00}       & 0.00                         & \multicolumn{1}{c|}{0.00}       & 0.00                         \\ \cline{2-10} 
				& 0.2                                              & \multicolumn{1}{c|}{0.00}       & 0.16                & \multicolumn{1}{c|}{0.00}       & 0.05       & \multicolumn{1}{c|}{0.00}       & 0.05                         & \multicolumn{1}{c|}{0.00}       & 0.05                         \\ \cline{2-10} 
				& 0.4                                              & \multicolumn{1}{c|}{0.26}       & 5.64                & \multicolumn{1}{c|}{0.08}       & 2.10       & \multicolumn{1}{c|}{0.08}       & 1.97                         & \multicolumn{1}{c|}{0.07}       & 1.88                         \\ \cline{2-10} 
				& 0.6                                              & \multicolumn{1}{c|}{3.28}       & 20.80               & \multicolumn{1}{c|}{1.16}       & 8.47       & \multicolumn{1}{c|}{1.09}       & 7.99                         & \multicolumn{1}{c|}{1.04}       & 7.66                         \\ \cline{2-10} 
				& 0.7                                              & \multicolumn{1}{c|}{7.54}       & 29.71               & \multicolumn{1}{c|}{2.77}       & 12.58      & \multicolumn{1}{c|}{2.61}       & 11.89                        & \multicolumn{1}{c|}{2.49}       & 11.41                        \\ \cline{2-10} 
				& 0.8                                              & \multicolumn{1}{c|}{13.70}      & 38.28               & \multicolumn{1}{c|}{5.25}       & 16.82      & \multicolumn{1}{c|}{4.94}       & 15.91                        & \multicolumn{1}{c|}{4.73}       & 15.29                        \\ \hline 
			\end{tabular}
	\end{center}
\end{table}
\newpage
\begin{table}[h!]
	\caption{Percentage risk improvement of $d_{11}^E$ and $d_{BZ}^E$ with respect to entropy loss function.}\label{Tab2}
\begin{center}
	\begin{tabular}{|c|c|cc|cc|cc|cc|}
		\hline
		\multirow{2}{*}{$(n_1,n_2)$} & \multirow{2}{*}{$\frac{\sigma_1^2}{\sigma_2^2}$} & \multicolumn{2}{c|}{$\nu=5$}                          & \multicolumn{2}{c|}{$\nu=8$}                 & \multicolumn{2}{c|}{$\nu=10$}                                  & \multicolumn{2}{c|}{$\nu=15$}                                  \\ \cline{3-10} 
		&                                                  & \multicolumn{1}{c|}{$d_{11}^E$} & \textbf{$d_{BZ}^E$} & \multicolumn{1}{c|}{$d_{11}^E$} & $d_{BZ}^E$ & \multicolumn{1}{c|}{$d_{11}^E$} & $d_{BZ}^E$  & \multicolumn{1}{c|}{$d_{11}^E$} & $d_{BZ}^E$ \\ \hline
		\multirow{6}{*}{$(5,7)$}     & 0.1                                              & \multicolumn{1}{c|}{0.04}       & 2.67                & \multicolumn{1}{c|}{0.02}       & 1.18       & \multicolumn{1}{c|}{0.01}       & 1.14                         & \multicolumn{1}{c|}{0.01}       & 1.09                         \\ \cline{2-10} 
		& 0.2                                              & \multicolumn{1}{c|}{0.52}       & 10.37               & \multicolumn{1}{c|}{0.22}       & 4.74       & \multicolumn{1}{c|}{0.21}       & 4.59                         & \multicolumn{1}{c|}{0.20}       & 4.41                         \\ \cline{2-10} 
		& 0.4                                              & \multicolumn{1}{c|}{5.08}       & 27.54               & \multicolumn{1}{c|}{2.22}       & 13.32      & \multicolumn{1}{c|}{2.15}       & 12.91                        & \multicolumn{1}{c|}{2.06}       & 12.42                        \\ \cline{2-10} 
		& 0.6                                              & \multicolumn{1}{c|}{14.28}      & 41.10               & \multicolumn{1}{c|}{6.48}       & 20.75      & \multicolumn{1}{c|}{6.27}       & 20.13                        & \multicolumn{1}{c|}{6.02}       & 19.38                        \\ \cline{2-10} 
		& 0.7                                              & \multicolumn{1}{c|}{19.66}      & 46.45               & \multicolumn{1}{c|}{9.07}       & 23.88      & \multicolumn{1}{c|}{8.79}       & 23.17                        & \multicolumn{1}{c|}{8.45}       & 22.32                        \\ \cline{2-10} 
		& 0.8                                              & \multicolumn{1}{c|}{25.13}      & 51.01               & \multicolumn{1}{c|}{11.80}      & 26.65      & \multicolumn{1}{c|}{11.43}      & 25.87                        & \multicolumn{1}{c|}{10.99}      & 24.92                        \\ \hline
		\multirow{6}{*}{$(10,8)$}    & 0.1                                              & \multicolumn{1}{c|}{0.00}       & 0.17                & \multicolumn{1}{c|}{0.00}       & 0.11       & \multicolumn{1}{c|}{0.00}       & 0.10                         & \multicolumn{1}{c|}{0.00}       & 0.09                         \\ \cline{2-10} 
		& 0.2                                              & \multicolumn{1}{c|}{0.06}       & 1.64                & \multicolumn{1}{c|}{0.04}       & 1.09       & \multicolumn{1}{c|}{0.03}       & 1.03                         & \multicolumn{1}{c|}{0.03}       & 0.92                         \\ \cline{2-10} 
		& 0.4                                              & \multicolumn{1}{c|}{1.10}       & 8.68                & \multicolumn{1}{c|}{0.73}       & 5.92       & \multicolumn{1}{c|}{0.68}       & 5.57                         & \multicolumn{1}{c|}{0.61}       & 4.98                         \\ \cline{2-10} 
		& 0.6                                              & \multicolumn{1}{c|}{4.53}       & 17.08               & \multicolumn{1}{c|}{3.05}       & 11.83      & \multicolumn{1}{c|}{2.86}       & 11.14                        & \multicolumn{1}{c|}{2.55}       & 9.99                         \\ \cline{2-10} 
		& 0.7                                              & \multicolumn{1}{c|}{7.11}       & 20.97               & \multicolumn{1}{c|}{4.81}       & 14.61      & \multicolumn{1}{c|}{4.52}       & 13.78                        & \multicolumn{1}{c|}{4.04}       & 12.37                        \\ \cline{2-10} 
		& 0.8                                              & \multicolumn{1}{c|}{10.08}      & 24.54               & \multicolumn{1}{c|}{6.87}       & 17.19      & \multicolumn{1}{c|}{6.46}       & 16.21                        & \multicolumn{1}{c|}{5.77}       & 14.57                        \\ \hline
		\multirow{6}{*}{$(12,15)$}   & 0.1                                              & \multicolumn{1}{c|}{0.00}       & 0.01                & \multicolumn{1}{c|}{0.00}       & 0.01       & \multicolumn{1}{c|}{0.00}       & 0.01                         & \multicolumn{1}{c|}{0.00}       & 0.01                         \\ \cline{2-10} 
		& 0.2                                              & \multicolumn{1}{c|}{0.00}       & 0.41                & \multicolumn{1}{c|}{0.00}       & 0.26       & \multicolumn{1}{c|}{0.00}       & 0.24                         & \multicolumn{1}{c|}{0.00}       & 0.22                         \\ \cline{2-10} 
		& 0.4                                              & \multicolumn{1}{c|}{0.28}       & 5.27                & \multicolumn{1}{c|}{0.18}       & 3.48       & \multicolumn{1}{c|}{0.16}       & 3.21                         & \multicolumn{1}{c|}{0.15}       & 2.95                         \\ \cline{2-10} 
		& 0.6                                              & \multicolumn{1}{c|}{2.59}       & 14.26               & \multicolumn{1}{c|}{1.69}       & 9.61       & \multicolumn{1}{c|}{1.55}       & 8.88                         & \multicolumn{1}{c|}{1.42}       & 8.18                         \\ \cline{2-10} 
		& 0.7                                              & \multicolumn{1}{c|}{5.14}       & 19.11               & \multicolumn{1}{c|}{3.38}       & 12.98      & \multicolumn{1}{c|}{3.11}       & 12.01                        & \multicolumn{1}{c|}{2.85}       & 11.07                        \\ \cline{2-10} 
		& 0.8                                              & \multicolumn{1}{c|}{8.54}       & 23.79               & \multicolumn{1}{c|}{5.65}       & 16.27      & \multicolumn{1}{c|}{5.21}       & 15.08                        & \multicolumn{1}{c|}{4.79}       & 13.91                        \\ \hline
		\multirow{6}{*}{$(13,18)$}   & 0.1                                              & \multicolumn{1}{c|}{0.00}       & 0.01                & \multicolumn{1}{c|}{0.00}       & 0.00       & \multicolumn{1}{c|}{0.00}       & 0.00                         & \multicolumn{1}{c|}{0.00}       & 0.00                         \\ \cline{2-10} 
		& 0.2                                              & \multicolumn{1}{c|}{0.00}       & 0.33                & \multicolumn{1}{c|}{0.00}       & 0.17       & \multicolumn{1}{c|}{0.00}       & 0.15                         & \multicolumn{1}{c|}{0.00}       & 0.13                         \\ \cline{2-10} 
		& 0.4                                              & \multicolumn{1}{c|}{0.25}       & 5.81                & \multicolumn{1}{c|}{0.12}       & 3.08       & \multicolumn{1}{c|}{0.11}       & 2.72                         & \multicolumn{1}{c|}{0.09}       & 2.36                         \\ \cline{2-10} 
		& 0.6                                              & \multicolumn{1}{c|}{2.88}       & 17.39               & \multicolumn{1}{c|}{1.48}       & 9.56       & \multicolumn{1}{c|}{1.30}       & 8.50                         & \multicolumn{1}{c|}{1.13}       & 7.39                         \\ \cline{2-10} 
		& 0.7                                              & \multicolumn{1}{c|}{6.14}       & 23.85               & \multicolumn{1}{c|}{3.21}       & 13.32      & \multicolumn{1}{c|}{2.84}       & 11.87                        & \multicolumn{1}{c|}{2.45}       & 10.34                        \\ \cline{2-10} 
		& 0.8                                              & \multicolumn{1}{c|}{10.70}      & 30.12               & \multicolumn{1}{c|}{5.69}       & 17.08      & \multicolumn{1}{c|}{5.04}       & 15.24                        & \multicolumn{1}{c|}{4.37}       & 13.30                        \\ \hline
		\multirow{6}{*}{$(20,17)$}   & 0.1                                              & \multicolumn{1}{c|}{0.00}       & 0.00                & \multicolumn{1}{c|}{0.00}       & 0.00       & \multicolumn{1}{c|}{0.00}       & 0.00                         & \multicolumn{1}{c|}{0.00}       & 0.00                         \\ \cline{2-10} 
		& 0.2                                              & \multicolumn{1}{c|}{0.00}       & 0.05                & \multicolumn{1}{c|}{0.00}       & 0.03       & \multicolumn{1}{c|}{0.00}       & 0.03                         & \multicolumn{1}{c|}{0.00}       & 0.03                         \\ \cline{2-10} 
		& 0.4                                              & \multicolumn{1}{c|}{0.07}       & 1.82                & \multicolumn{1}{c|}{0.05}       & 1.30       & \multicolumn{1}{c|}{0.04}       & 1.20                         & \multicolumn{1}{c|}{0.04}       & 1.11                         \\ \cline{2-10} 
		& 0.6                                              & \multicolumn{1}{c|}{1.03}       & 7.21                & \multicolumn{1}{c|}{0.73}       & 5.22       & \multicolumn{1}{c|}{0.67}       & 4.84                         & \multicolumn{1}{c|}{0.62}       & 4.48                         \\ \cline{2-10} 
		& 0.7                                              & \multicolumn{1}{c|}{2.43}       & 10.70               & \multicolumn{1}{c|}{1.73}       & 7.78       & \multicolumn{1}{c|}{1.60}       & 7.23                         & \multicolumn{1}{c|}{1.48}       & 6.69                         \\ \cline{2-10} 
		& 0.8                                              & \multicolumn{1}{c|}{4.57}       & 14.31               & \multicolumn{1}{c|}{3.28}       & 10.45      & \multicolumn{1}{c|}{3.04}       & 9.72                         & \multicolumn{1}{c|}{2.81}       & 9.00                         \\ \hline
	\end{tabular}
\end{center}
\end{table}
		\newpage
	\newpage
\section{Conclusions}
	In this manuscript, we studied component-wise estimation of the ordered variances of two scale mixture of normal distributions. In this study, we focused on developing improved estimators for $\sigma_1^2$ over the BAEE. Moreover, we have derived a class of improved estimators of BAEE of $\sigma_1^2$. This study was carried out specifically under two loss functions, namely  (i) the squared error loss function, (ii) the entropy loss function. We have also derived similar improvement results for $\sigma_2^2$. Furthermore, as an application improvement results are obtained for estimating ordered scale parameters of a scale mixture of multivariate normal distributions. Finally, we have conducted the simulation study and tabulated the percentage risk improvement with respect to BAEE, which validates our findings numerically. 
\bibliography{bib_normal}
\end{document}